\documentclass[preprint,12pt]{elsarticle}

\usepackage{graphicx}
\usepackage{epstopdf}
\usepackage{amsmath,amssymb}
\usepackage{rotating}
\usepackage{sidecap}
\usepackage{natbib}
\usepackage{color}

\numberwithin{equation}{section}

\begin{document}

\begin{frontmatter}

\title{Model of the Human Sleep Wake System}

\author{Lisa Rogers\corref{lisa}\fnref{lisa2}}
\ead{lrogers@cims.nyu.edu}

 \address{Courant Institute of Mathematical Sciences, New York University, NY, NY}
 \author{Mark Holmes }
 \ead{holmes@rpi.edu}
 \address{Rensselaer Polytechnic Institute, Troy, NY}

\begin{abstract}
{A model and analysis of the human sleep/wake system is presented. The model is derived using the known neuronal groups, and their various projections, involved with sleep and wake.  Inherent in the derivation is the existence of a slow time scale associated with homeostatic regulation, and a faster time scale associated with the dynamics within the sleep phase.  A significant feature of the model is that it does not contain a periodic forcing term, common in other models, reflecting the fact that sleep/wake is not dependent upon a diurnal stimulus.  Once derived,  the model is analyzed using  a linearized stability analysis. We then use experimental data from normal sleep-wake systems and orexin knockout systems to verify the physiological validity of the equations.}
\end{abstract}
\end{frontmatter}

\section{Introduction}
Most living organisms exhibit daily biological rhythms. One of the most curious and difficult to understand is the sleep/wake cycle. Even though we are very familiar with sleep/wake, the reasons for sleep are not well understood. What is intriguing about the sleep-wake cycle is the complex interplay between the arousal and sleep states. Part of this is due to the cyclic alteration of the rapid eye movement (REM) and non-REM (NREM) stages of sleep.  The distribution of REM and NREM varies dramatically across the animal kingdom, to the point that a clear pattern has not emerged.  One consequence of this is that an animal model for human sleep has not been found \cite{Zimmer}. Another complication is that sleep/wake is sometimes thought to be synonymous with the 24 hour circadian clock.  This is not correct.  Although sleep/wake is coupled to the circadian clock, it is regulated by a neural system distinct from the circadian timing mechanism \cite{Semba}. Sleep/wake are also heavily influenced by two other factors.  One is \textit{homeostatic}, where the need for sleep increases with the length of the wake period.  The second is \textit{allostatic}, where the sleep/wake cycle is modified due to external behavioral events \cite{saper}. The circadian process is connected to the sleep-wake cycle, while a homeostatic process is inherent to the sleep-wake cycle. To be specific, it is now known that the sleep-wake cycle is independent of the diurnal stimulus, and is instead driven by homeostatic forces \cite{scammell,Jones,borbely}. Therefore, the inherent periodicity of the sleep/wake cycle is not due to a periodic forcing mechanism associated with the diurnal cycle but is an inherent characteristic of the  system.  The objective of this study is to derive, and then analyze, a model for the sleep/wake cycle that accounts for the circadian and homeostatic influences of sleep/wake, but not necessarily the allostatic affects.  This will include both REM and NREM sleep, as well as the waking state and observations about the associated brainwave activity of all states. 

\begin{figure}[h]
	\begin{center}
		\includegraphics[width=5in]{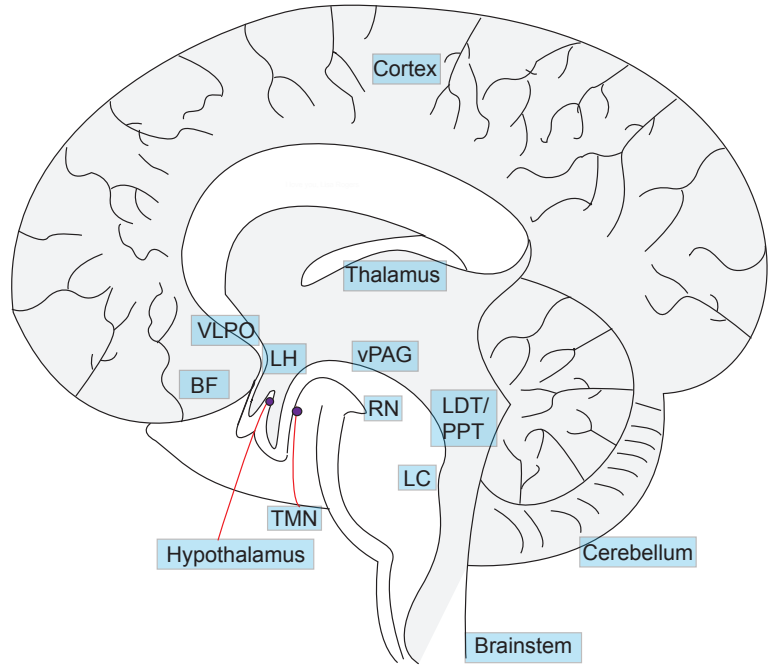}
	\caption{Approximate locations of areas pertaining to Sleep-Wake promoting neuronal groups.}
	\label{fig:fig4}
 \end{center}
\end{figure}

\subsection{Biological Background}
A basic assumption is that the brain has two states, sleep and wake, with sharp transitions between them.  There are multiple cell groups in the brain contributing to these states, but the monoaminergic nuclei (MN) have been identified as primarily responsible for promoting the awake state \cite{saper}.  Sleep is primarily promoted by the ventrolateral preoptic nucleus (VLPO), which is a group of cells associated with generating NREM and REM. A consequence of this distinct structural organization is that although the states are coupled, they have different regulatory mechanisms in the brain.  

\begin{figure}[h!]
\begin{center}
\includegraphics[width=4.5in]{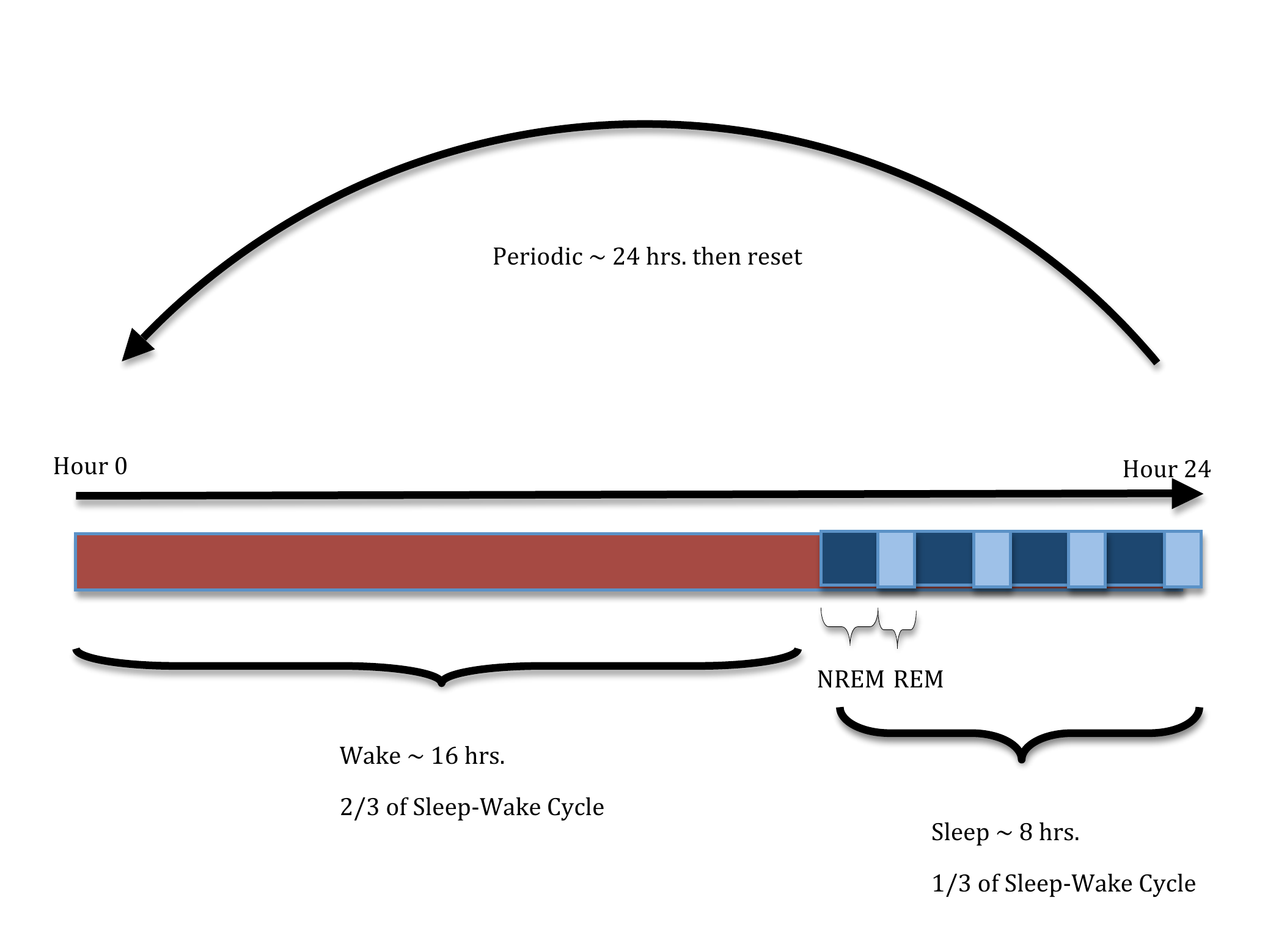}
\caption{The basic partitions of the states of sleep and wake for humans.}
\label{fig:fig1}
\end{center}
\end{figure}

\subsubsection{REM/NREM}

The REM and NREM phases of sleep have a well defined pattern (refer to Figure \ref{fig:fig1}).  Regardless of the conditions of sleep onset, a person always begins in an NREM period of the cycle.  After approximately 80-85 minutes the brain switches to the REM stage, and remains in this state for approximately 5-10 minutes before switching back to NREM.  This 90 minute pattern continues through the night with the REM-on intervals becoming progressing longer, and the NREM-on intervals shorter.  The last REM-on period is about 25 minutes.  There is no probability of waking naturally during an NREM interval, as the NREM component occurs during the periods of slow wave sleep. It is postulated that this is due to NREM being the most mentally restorative of the two processes, requiring little brain function and unconsciousness. Natural wake up occurs towards the end of the last NREM period. \cite{SRS} \\

There are three different ways of triggering sleep onset: exciting sleep active neurons, inhibiting wake active neurons or a combination of the two. The same idea applies to the triggering of wake onset.  Consequently, to investigate the source of sleep/wake brain function, it is necessary to examine the neurotransmitter activity of the system. 

\subsubsection{Location Specific Neurotransmitter Activity}
\begin{figure}[h]
	\begin{center}
		\includegraphics[width=5in]{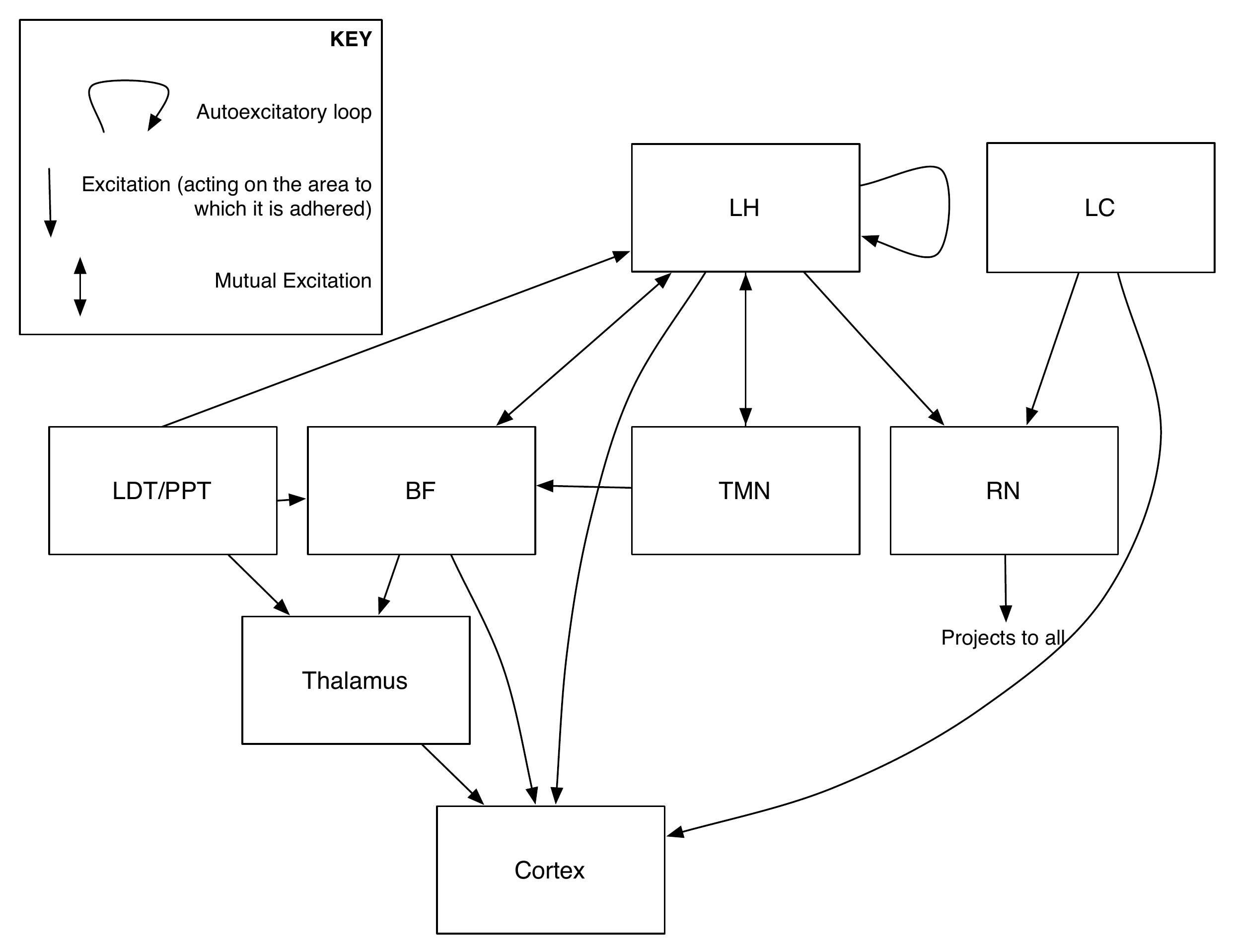}
	\caption{Excitatory projections of neurotransmitter flow in the sleep/wake generating network.}
	\label{fig:fig2}
 \end{center}
\end{figure}

\begin{figure}[h!]
	\begin{center}
		\includegraphics[width=5in]{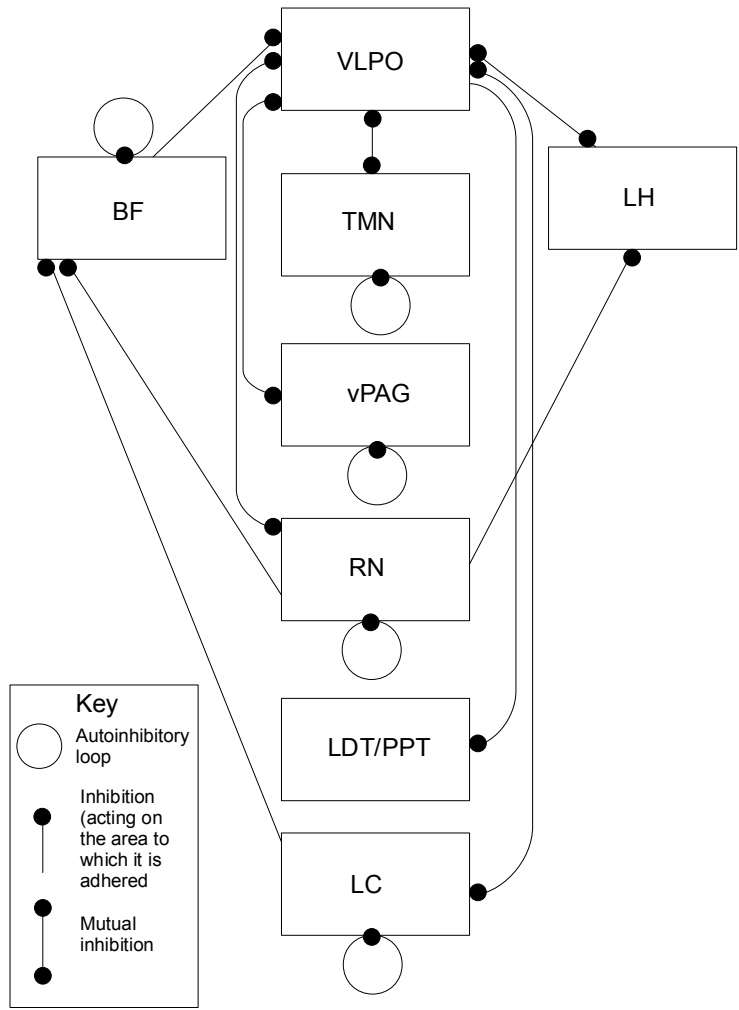}
	\caption{Inhibitory projections of neurotransmitter flow in the sleep/wake generating network}
	\label{fig:fig3}
 \end{center}
\end{figure}

\begin{figure}[h]
	\begin{center}
		\includegraphics[width=5in]{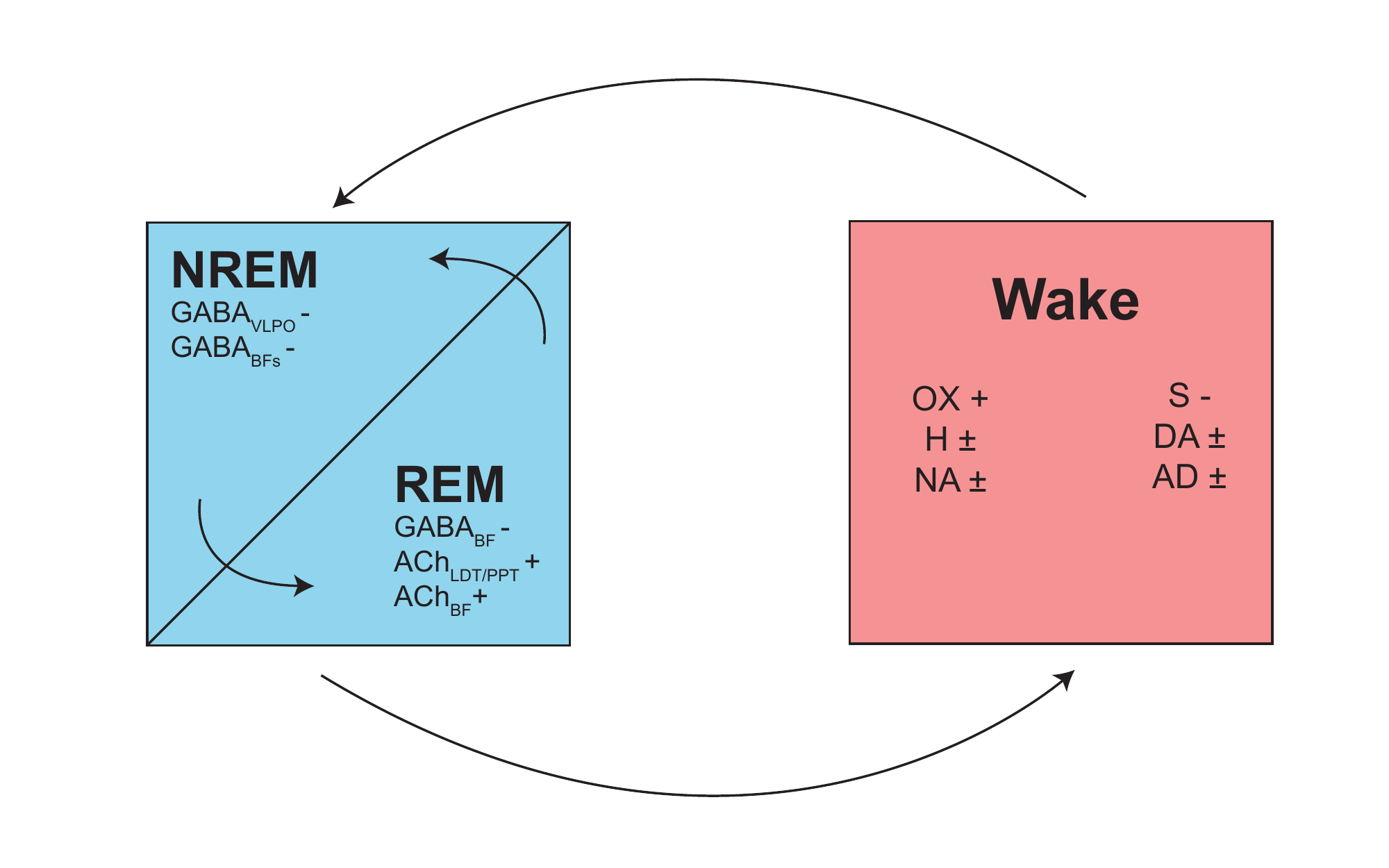}
	\caption{The cycle of NREM Sleep, REM Sleep and Wakefulness with the neurotransmitters characteristic of their stages.}
	\label{fig:fig5}
 \end{center}
\end{figure}

The chemical activity of the brain involves the secretion of hormones and neurotransmitters from specific sections of the brain. The sleep-wake cycle is initiated and maintained by multiple distinct neuronal regions in the brainstem. These regions are located in the hypothalamus and are coupled to each other via electrical and neurochemical signaling. Figures \ref{fig:fig2} and \ref{fig:fig3} provide a schematic of location specific neurotransmitter signal flow in terms of excitatory and inhibitory mechanisms. Figure \ref{fig:fig4} shows the approximate locations of sleep and wake promoting areas in the human brain. Figure \ref{fig:fig5} provides a simple depiction of the switching between wake and sleep states and well as the switch between REM and NREM states, and the set of neurotransmitters associated with each, as well as their excitatory (+) or inhibitory (-) effects. The major neurotransmitter players are GABA (gamma-aminobutyric acid), orexin (hypocretin), histamine, acetylcholine, noradrenaline (norepinephrine), serotonin and dopamine. Each of these transmitters are located throughout the following neuronal groups: basal forebrain (BF), raphe nuclei (RN), laterodorsal tegmentum (LDT), pedunculopontine tegmentum (PPT), locus coeruleus (LC), lateral hypothalamus (LH), tuberomammillary nucleus (TMN), ventral periaqueductal gray region (vPAG), and ventrolateral preoptic nucleus (VLPO). Adenosine has been proposed to be the homeostatic accumulator of the need for sleep \cite{saper}. Acetylcholine effects muscle action, learning and memory. Norepinephrine effects one's control of alertness and wakefulness. Dopamine effects movement, focus and learning. Serotonin regulates mood, appetite, sleep, arousal and pain threshold. GABA contributes to sleep control and appetite control.  Also note that the sections of the brain we will be focusing on are the cerebral cortex, hypothalamus and thalamus, which play key roles in homeostatic maintenance. In this way they are the most relevant to the control of the human sleep/wake system. \\

\begin{figure}[h]
	\begin{center}
		\includegraphics[width=5in]{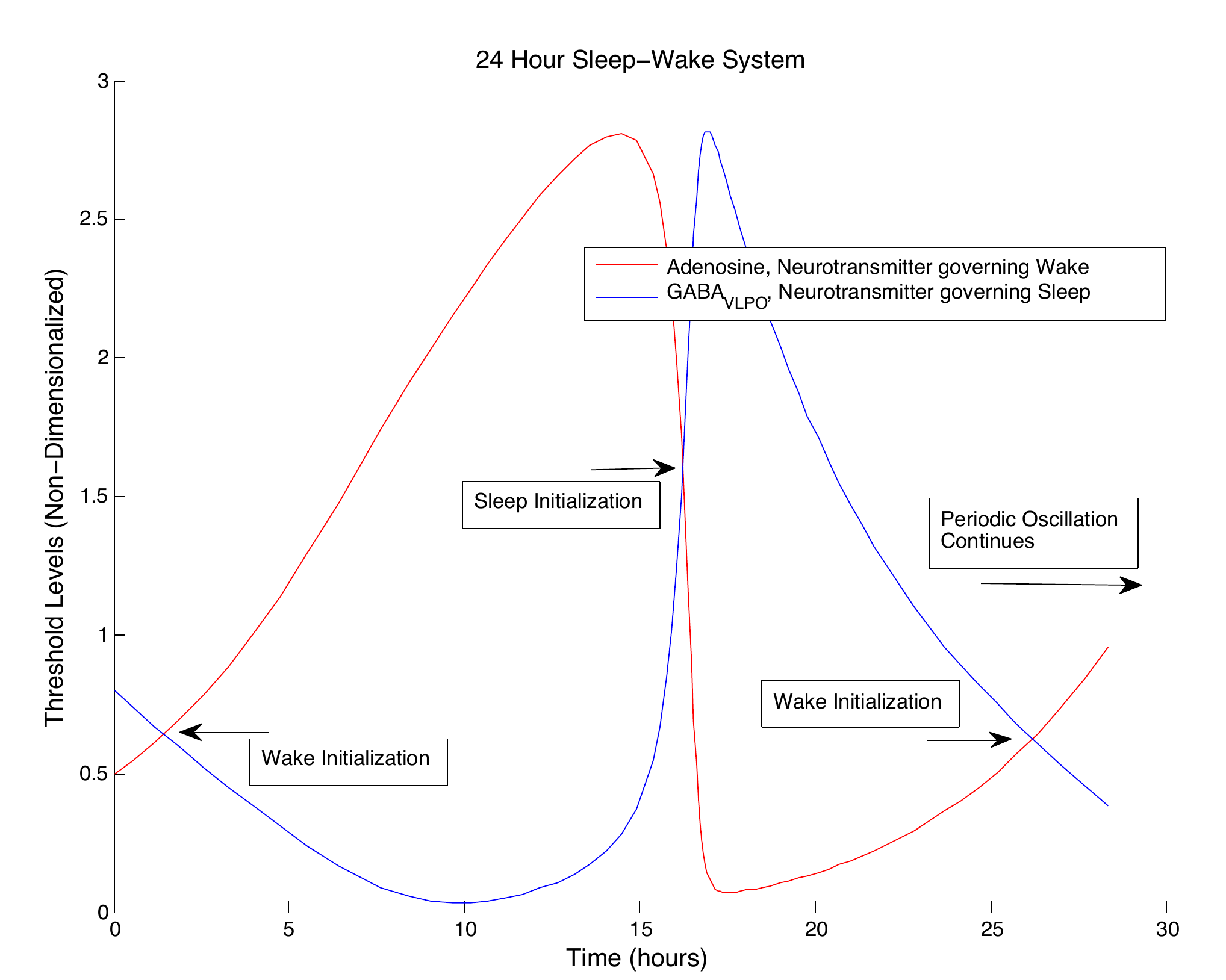}
	\caption{One full period of the sleep-wake cycle model.}
	\label{fig:fig6}
 \end{center}
\end{figure}

Figure \ref{fig:fig6} shows one full period of the sleep-wake cycle model in terms of adenosine and GABA concentration levels. The point at which the curves cross indicate wake or sleep initializations, depending on whether one of the variables is increasing or decreasing. Wake Initialization occurs when $GABA$ is decreasing and $AD$ is increasing. Sleep initialization occurs when $AD$ is decreasing and $GABA$ is increasing. 

\newpage
\begin{sidewaystable}[h]
\centering
\vspace{0.25in}
\caption{Neurotransmitter Details}
\bigskip
\begin{tabular}{| p{0.7in} | p{0.34in}| p{0.7in} | p{1.15in} | p{1.15in} | p{0.8in} | p{0.7in} | p{0.7in} |}
\hline
Variable\footnote{Citations for the information found in this table are the following: \cite{adenosine1,Principles,scammell,adenosine_effects,chou_lu,arousal_systems,serotonin,Jones2,Dzirasa,Jones,Baumann,Arrigoni,adenosine2,Wu,Alreja,Khateb,Khateb2,John}} &  Region & Neurotrans -mitter & Receives Projections From: & Projects to: & Wake & NREM & REM\\
\hline
$GABA_{VLPO}$ & VLPO & GABA & LC, vPAG, RN, BF & LH, LC, LDT/PPT, vPAG, RN, BF & Inactive & High firing rate, high level. & Low firing rate, low level.  \\
\hline
$GABA_{BFs}$ & BF & GABA & BF, RN, LC, LH & BF, VLPO & Inactive & High firing rate, high level. & Inactive\\
\hline
$GABA_{BFw}$ & BF & GABA & BF, RN, LC, LH & BF, VLPO & Active & Inactive & Active \\
\hline
$OX$ & LH & Orexin & LH, VLPO, RN, TMN & LH, TMN, RN, BF & High firing rate, high level. & Low firing rate, low level. & Inactive \\
\hline
$H$ & TMN & Histamine & TMN, LH, VLPO & TMN, LH & Tonic, low firing rate, high level. & Low firing rate, low level. & Inactive \\
\hline
$ACh_{BF}$ & BF & Acetylcholine & BF, RN, LC, LH & BF, VLPO & High firing rate, high level. & Low firing rate, low level. & High firing rate, high level. \\
\hline
$ACh_{LDT/PPT}$ & LDT/ PPT & Acetylcholine & BF, VLPO & VLPO & Fast cortical rhythms & Inactive & Fast cortical rhythms. \\
\hline
 $NA$ & LC & Noradrenaline & LC, VLPO & LC, BF, RN, VLPO & Tonic, low firing rate, high level. & Low firing rate, low level. & Inactive \\
  \hline
 $S$ & RN & Serotonin & RN, LH, LC, VLPO & RN, VLPO, LH, BF & Tonic, low firing rate, high level. & Low firing rate, low level. & Inactive \\
\hline
$DA$ & vPAG & Dopamine & vPAG, VLPO & vPAG, VLPO & No rate change. Elevated levels. & No rate change. & No rate change. \\
\hline
$AD$ & All Cells & Adenosine & Not Applicable & Not Applicable & Levels increase with prolonged wake in the BF. & Levels fall during recovery sleep. & Levels fall during recovery sleep. \\
\hline
\end{tabular}
\end{sidewaystable}

\begin{figure}[t]
	\begin{center}
		\includegraphics[width=4in]{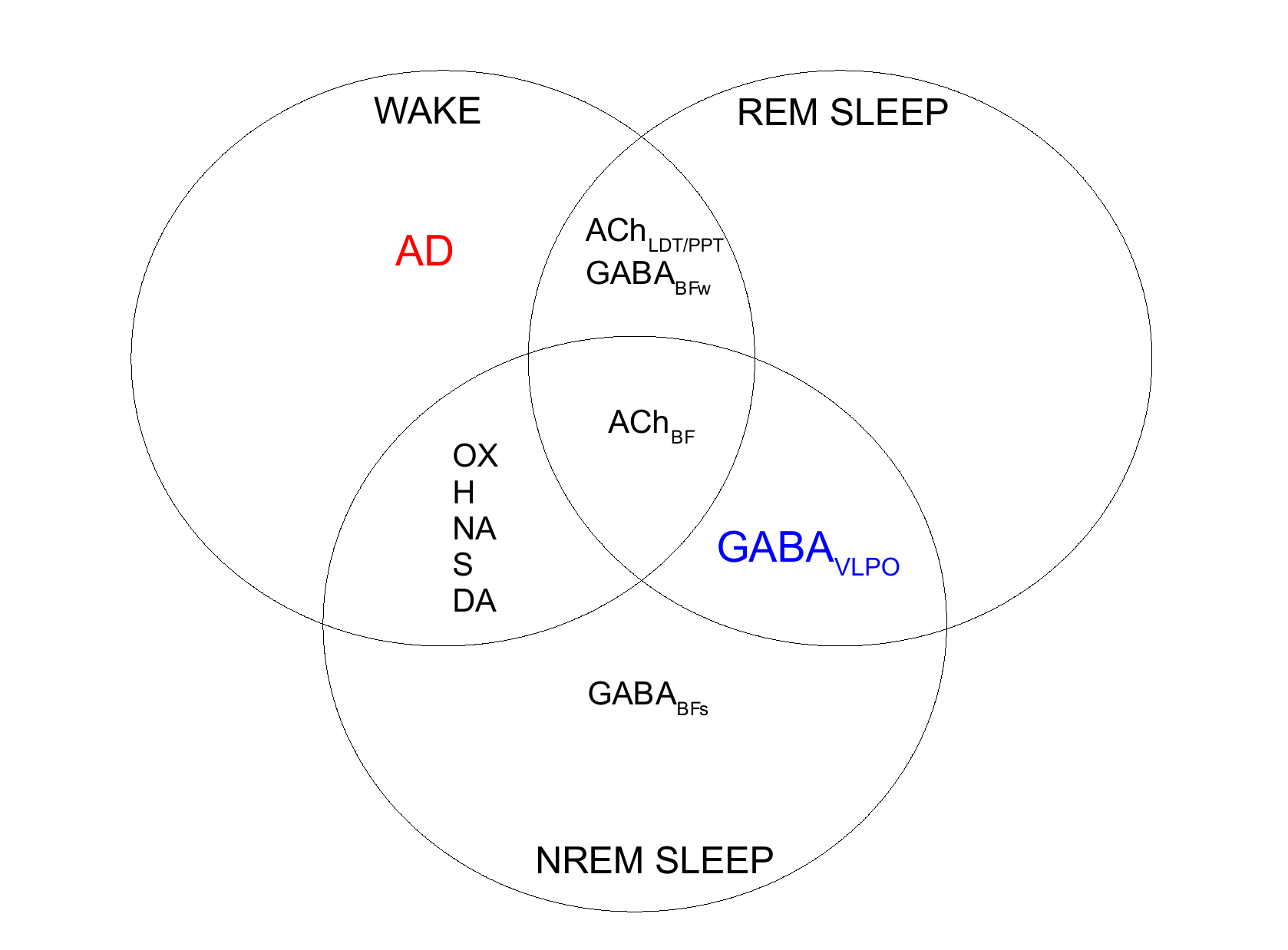}
	\caption{Components of the human sleep/wake cycle as well as their associate neurotransmitters}
	\label{fig:fig7}
 \end{center}
\end{figure}

Adenosine, which is a neuromodulator, plays a major role in the sleep-wake process. It is a byproduct of adenosine-triphosphate (ATP), the main ``energy currency" of the brain, whenever neurotransmitters are fired, adenosine is produced. Adenosine is also continuously converted back to ATP via the metabolic processes that use it as an energy source \cite{metabo}. The adenosine concentration that we are most interested in is that of the basal forebrain. While inhibiting the expression of all wake active neurotransmitters and disinhibiting the expression of those found in the VLPO \cite{adenosine_effects}, adenosine levels in the basal forebrain increase during prolonged wake and fall during recovery sleep, providing an intrinsic periodic force associated with waking \cite{adenosine1,adenosine2}. The linear relationship between between adenosine accumulation and sleep propensity provides support for using adenosine as representative of sleep propensity. GABA in the VLPO can be viewed as the other major player in the sleep-wake system. $GABA_{A}$ receptors in the brain mediate the \textit{tonic} (or slow) conductance \cite{Winsky}. This slower type of signaling plays a key role in cell excitability and sleep homeostasis \cite{Winsky}. Neuron populations in the VLPO are inactive during the vigilance state, and have a high firing rate and level during NREM, especially during delta wave (slow wave) sleep \cite{Winsky}.  Being the major inhibitory neurotransmitter in the central nervous system, GABA is a key player in many physical systems. 

Referring to the Table 1, we can see the detailed activity associated with each neurotransmitter and brain area, allowing for the development of a system of variables. As a simple example, $OX_{LH}$ is the variable associated with the chemical concentration of orexin in the lateral hypothalamus. Since we only discuss orexin in one location, the subscript LH can be dropped. Orexin expression is promoted by acetylcholine in the LDT/PPT and noradrenaline in the LC. Orexin expression is inhibited by adenosine, serotonin, GABA and noradrenaline. Orexin promotes the release of noradrenaline, serotonin, histamine, GABA in the basal forebrain, acetylcholine, and dopamine, as well as promoting its own release (autoexcitation). These dynamics indicate the need for a linear combination of acetylcholine, noradrenaline, serotonin, dopamine, histamine, GABA, adenosine, and orexin to express the rate of change of OX concentration.\\

\section{Materials and Methods}
In order to test the functionality of the model, we use numerical methods to determine how the system responds to a variety of perturbations.  In the calculations to follow, MATLAB 2009a was used to run all the programs written for the Human Sleep/Wake system described here. However, a fourth order Runge-Kutta method was used for solving the system of equations, instead of using one of the MATLAB ode solvers because of the unique perturbations used in the sleep experiments described.  

\section{Theory and Calculation}

\subsection{Previous Work}
A good survey of the early sleep/wake models can be found in \cite{strogatz2}.  We will only comment here on those most relevant to our development. \\

In a recent paper by Phillips and Robinson \cite{Phillips}, a model was proposed in order to improve quantitative understanding of the human sleep/wake cycle, which is similar to the intentions of this study. They use the following three neuronal populations to model the human ascending arousal system: monoaminergic (MA) nuclei, acetylcholine related (ACh) nuclei, and populations in the VLPO. The VLPO is then split into circadian and homeostatic drives, which are represented by periodic forcing functions. Each population has its dynamics modeled by a sum of mean cell body potentials and mean neuronal firing rates subjected to the periodic forcing of the VLPO. While their model consistently predicts sleep/wake cycle dynamics observed in previous studies, it is important to note that the homeostatic and circadian drives are represented as periodic forcing rather than dynamics inherent to the physical system. \\

\cite{Cecilia} developed a model of the network of neurons controlling behavioral state transitions in mice, while also identifying special features of human and mouse sleep. It is similar to the Phillips and Robinson model (as well as the model being developed here) in that it defines three main processes being modeled: Sleep activity, Wake activity, and REM activity. Each population is defined by a set of three equations, determined by a wide range in timescales. It incorporates homeostatic REM and NREM sleep drive and uses known sleep active neurotransmitters to determine parameters. It then utilizes extensive experimental data gathered from mouse sleep/wake states to verify the accuracy of the model. A limitation of their approach concerns the separation of vigilance states, which are characterized in a non-quantitative fashion. Similarly, ``activity states" are used rather than specific neurotransmitter concentration and/or neuron firing data. \\ 

Taking a somewhat different viewpoint,  \cite{best} use the ``flip-flop" dynamics proposed in \cite{saper}  to construct a biologically-based mathematical model that accounts for several features of the human sleep/wake cycle.  This differs from our model in several ways, but two examples are that they assume, incorrectly, that the SCN projects directly to the VLPO, and an external periodic forcing mechanism is used to drive the model. \\

In another recent paper by \cite{DinizBehn2010}, a model of the rat sleep/wake regulatory network is formed utilizing similar techniques as those used to model a mouse sleep/wake regulatory network in her aforementioned work \cite{Cecilia}. In this model,  neurotransmitter concentrations  are coupled  to neuronal group firing rates to create a more complete view of the hypothalamic sleep-wake regulatory network. The model is composed of three parts: wake/sleep active neuron population firing rates, wake/sleep active neurotransmitter concentrations and a homeostatic regulator. All physical parameters were determined using extensive experimental data. This paper emphasizes the importance of utilizing neuron firing rates as well as neurochemical concentrations to form a physiologically accurate model of sleep/wake regulation. A downside of this model lies in the use of steady state saturation functions and Heaviside functions. Rather than using inherent neurotransmitter dynamics, the functions of neurotransmitter release are approximated using sigmoidal functions. The homeostatic drive is a physically unquantifiable force composed of Heaviside switching mechanisms. \\

\subsection{Model Formulation}

In this model for the sleep/wake cycle, we are considering three principal states.  They are the wake state ($W$), the NREM sleep state ($S$) and the REM sleep state ($R$). Table 1 provides a framework for the detailed interactions of all the neurotransmitters involved in the sleep-wake system.  Given the number of variables involved,  before deriving the equations it is worth describing the dynamics in more qualitative terms.  \\

The sleep-specific and wake specific neurotransmitters operate on a \textit{slow time scale}, measured in terms of hours, while the REM and NREM specific neurotransmitters operate on a \textit{fast time scale system}, measured in terms of minutes. The two neurochemicals associated with the slow time scale are adenosine in the BF and GABA in the VLPO. We use these two concentrations as variables indicating promotion of wake and sleep. High adenosine concentrations in the BF promote wakefulness, while high GABA concentrations and fast firing rates in the VLPO promote sleep.  Adenosine's slow build up and decay is associated with wakefulness, while GABA in the VLPO is associated with sleep. These observations are the basis for the hypothesis that adenosine is the principal component in determining wake, and $GABA_{VLPO}$ is the principal component in promoting sleep \cite{Winsky,adenosine1,adenosine2}. We also know from these references that AD disinhibits the VLPO during slow wave sleep, providing a direct connection between the building adenosine sleep pressure during waking and the gradual relief of this pressure during sleep. The other neurotransmitters are assumed to be associated with the fast time scale dynamics.  Adenosine and $GABA_{VLPO}$ fluctuate with a period of 24 hours, while the other neurotransmitters fluctuate with a period of several minutes. \\

An important inherent quality of the sleep-wake cycle is that it should maintain it's state transitions in the absence of external periodic forcing. Mathematically, this indicates the need for the system to have a \textit{stable (attracting) limit cycle}.   This solution should also be stable to certain perturbations, such as those that might occur from genetic mutations causing the production of too much or too little of a neurotransmitter, or lacking the proper neuron groups to produce a neurotransmitter at all. Perturbations may also occur via the removal of specific neuron groups, blockage of neurotransmitter bonding sites, and external forcing. It is well known from experimental data that perturbations to neurotransmitters don't alter the two-state nature of the sleep-wake system \cite{adenosine1,adenosine2,Arrigoni,OXKO,John,Lee,Khateb,Khateb2,arousal_systems,Dzirasa,Lu,Crochet,Gervasoni,Baumann}. As an example, cholinergic agonists applied to the LC in cats leads to increased REM sleep, but the sleep-wake cycle as a whole remains periodic \cite{Crochet}. Therefore one of the central goals will be to locate the area in the phase plane of $GABA_{VLPO}$ and $AD$ where a limit cycle exists and is stable.  \\

\subsubsection{Kinetics}
For each neurotransmitter involved in the sleep-wake system, we identify a location specific variable representing the effective neurotransmitter concentration. We assume first order kinetics for the neurotransmitters operating on the fast time scale. The fast time system may be formed from the observations compiled about the inhibitory and excitatory nature of the neurotransmitters in Table 1 combined with the dynamics arising in a synapse.  If $X_{i}$ is representative of a neurotransmitter in region $i$, we have 
\begin{equation}
\frac{dX_{i}}{dt} = {\rm Release(t)} - {\rm Uptake(t)} - {\rm Degradation(t)}. \label{eq:eq1}
\end{equation}
Note that inhibitory effects are considered pre-synaptic and are therefore controlling flow from the neuron into the synaptic gap. Despite the fact that inhibitory and excitatory effects can occur pre- or post-synaptically, for simplicity we will assume the blocking or binding process of inhibition as only controlling flow into the gap. The Release function contains rate parameters and neurotransmitter variables that excite and inhibit the expression of $X_{i}$. The Uptake function depends on $X_{i}$ and its associated vesicular uptake transporter, if existent. The Degradation function depends on $X_{i}$ and an associated rate of decay. The uptake transporters are modeled using known logistic dynamics from neurotransmitter uptake transporter microdialysis measurements \cite{neurotransport}. \\

As an example, the release of histamine is promoted by orexin, noradrenaline and acetylcholine in the BF and inhibited by GABA in the VLPO. Histamine is autoinhibitory, like most of the monoamines. The equation for histamine's rate of change, like all of the other neurochemicals involved, also requires an adenosine term. Histamine in the human sleep-wake system has no known uptake mechanism, and so only degrades as a method of removal. Therefore the equation for histamine is
\begin{equation}
 \frac{dH}{dt} = c_{i}OX + c_{i+1}NA + c_{i+2}ACh_{BF} + c_{i+3}AD - c_{i+4}GABA_{VLPO} - (c_{i+5} - c_{i+6})H \notag 
\end{equation}
where the $c_{i}$'s are rate coefficients. The rest of the equations are formulated in the same fashion.  

\bigskip\noindent
\textbf{Fast Time System}
\begin{align}
\frac{dGABA_{BFw}}{dt} &= c_{1}ACh_{BF} + c_{2}OX + c_{3}AD - c_{4}GABA_{BFw} - \notag \\
& \qquad c_{5}(GA1)(GABA_{BFw} ) \label{eq:eq2}\\
\frac{dGABA_{BFs}}{dt} &=  c_{6}AD + c_{7}NA - c_{8}GABA_{BFs} - c_{9}(GA1)(GABA_{BFs}) \label{eq:eq3}\\
\frac{dOX}{dt} &= c_{10}ACh_{LDT/PPT} + (c_{11} - c_{12})NA + (c_{13}-c_{14})OX + c_{15}AD - \notag \\
& \qquad c_{16}S - c_{17}GABA_{VLPO} + c_{18}GABA_{BFs} \label{eq:eq4}\\
\frac{dH}{dt} &= c_{19}OX + c_{20}NA + c_{21}ACh_{BF} + c_{22}AD - c_{23}GABA_{VLPO} - \notag \\
& \qquad (c_{24} - c_{25})H \label{eq:eq5} \\
\frac{dACh_{BF}}{dt} &= c_{26}OX + c_{27}H + c_{28}AD - c_{29}S + c_{30}GABA_{BFs} - c_{31}ACh_{BF} \label{eq:eq6} \\
\frac{dACh_{LDT/PPT}}{dt} &= c_{32}OX + c_{33}H + c_{34}AD + c_{35}NA - c_{36}S - c_{37}ACh_{LDT/PPT} \label{eq:eq7} \\
\frac{dNA}{dt} &= c_{38}OX + c_{39}ACh_{BF} + c_{40}AD - c_{41}GABA_{VLPO} - \notag \\
& \qquad (c_{42} - c_{43})NA - c_{44}(hNET)(NA) \label{eq:eq8} \\
\frac{dS}{dt} &= c_{45}OX + c_{46}H + c_{47}NA + c_{48}AD - c_{49}GABA_{VLPO} - \notag \\
& \qquad (c_{50}  - c_{51})S - c_{52}(hSERT)(S) \label{eq:eq9} \\
\frac{dDA}{dt} &= c_{53}OX + c_{54}S + c_{55}ACh_{LDT/PPT} + c_{56}NA + c_{57}AD - \notag \\
& \qquad c_{58}GABA_{VLPO} - c_{59}(hDAT)(DA) - c_{60}DA  \label{eq:eq10}
\end{align}

The uptake transporter values are contained in the uptake terms in equations (\ref{eq:eq2}), (\ref{eq:eq3}), (\ref{eq:eq8}), (\ref{eq:eq9}) and (\ref{eq:eq10}). Their values are found from microdialysis experiments with human cerebrospinal fluid (CSF), and are as follows \cite{neurotransport}: $GA1 = 1$,
$hNET = 0.457$, 
$hSERT = 0.463$, 
$hDAT = 1.22 $.\\

\subsubsection{Uptake Mechanisms}

The neurotransmitters involved in the sleep-wake system, as with any neuronal network, have natural rates of removal from the synaptic gap once released. These rates can either be contributed to by degradation, neurotransmitter uptake or a combination of both. All of the neurotransmitters involved degrade over time, but only some have neurotransmitter uptake mechanisms. We will not be considering any of the smaller scale chemical dynamics of the neurotransmitter uptake mechanisms. We care only whether a neurotransmitter is removed  by an uptake mechanism, and will not differentiate between plasma membrane-bound transportation or receptor uptake. The key take away from the uptake transportation data is that only GABA, noradrenaline, serotonin and dopamine have significant synaptic uptake transport mechanisms in humans. These mechanisms are characterized by the data given constants $GA1$, $hNET$, $hSERT$ and $hDAT$ for GABA, noradrenaline, serotonin and dopamine, respectively \cite{neurotransport}. These neurotransmitter transporters appear in the uptake term of the equations, always multiplied by the neurotransmitter it is dependent on to indicate a linear removal rate. As an example, $hSERT$ appears in equation 9 multiplied by $-c_{52}S$, so as to act as a removal mechanism for serotonin when it is present. However, $hSERT$ does not multiply the natural degradation term $-c_{50}S$.\\

\subsubsection{Nonlinear Dynamics}
\begin{figure}[t]
	\begin{center}
		\includegraphics[width=5in]{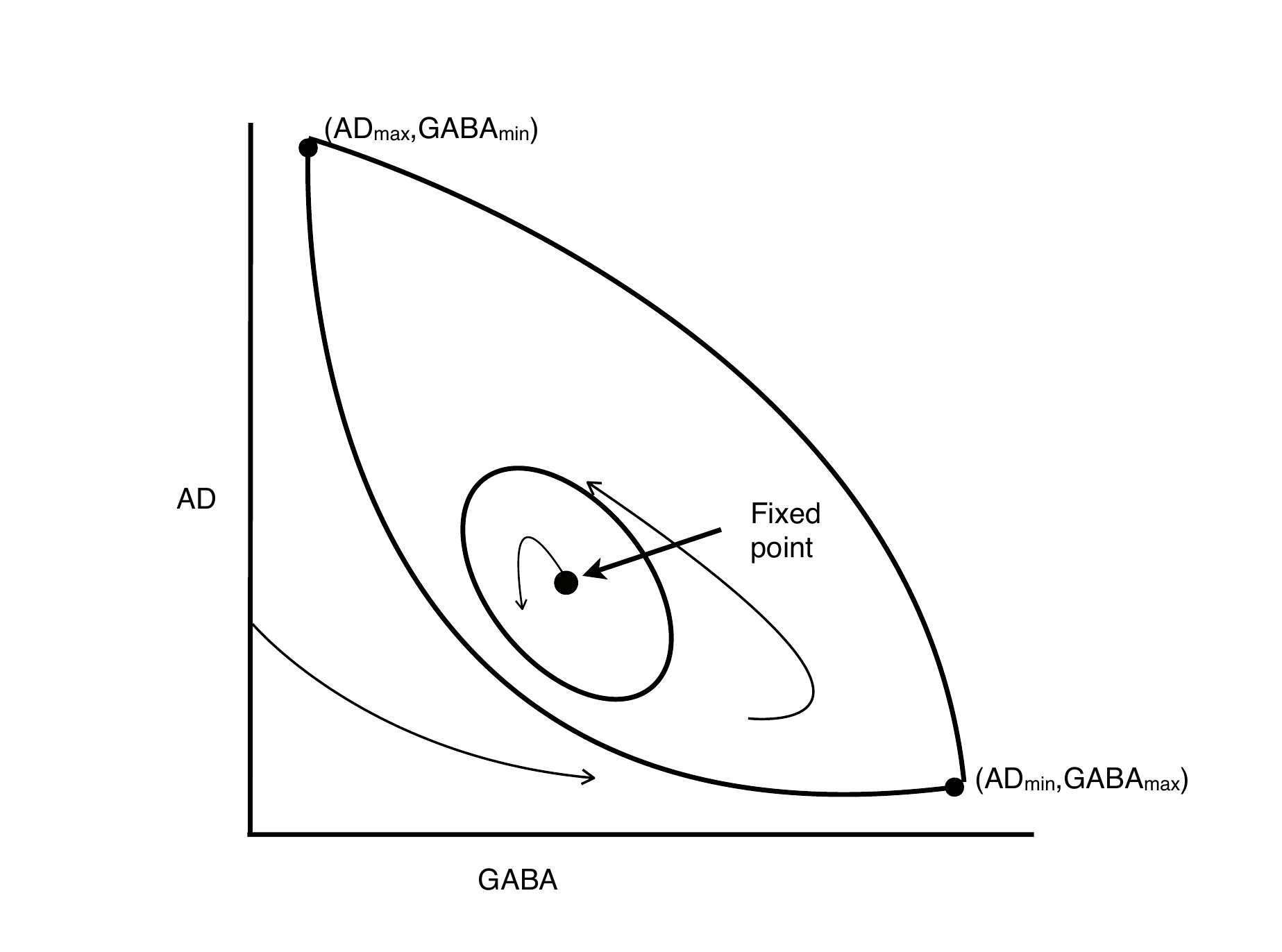}
	\caption{The ideal activity of the two neurotransmitter nullclines in the $GABA_{VLPO}-AD$ phase plane.}
	\label{fig:fig8}
 \end{center}
\end{figure}

There are two neurotransmitters, $AD$ and $GABA_{VLPO}$, operating on a slow time scale that define the sleep-wake cycle. The qualitative portrait in phase space is known from basic sleep-wake cycle activity. We can construct the desired closed orbit and concentric curves from the nullclines of the $AD$ and $GABA_{VLPO}$ equations. We identify the bounding region in the phase space activity by the box contained between $(AD, GABA_{VLPO})  = (AD_{max},GABA_{min})$ and $(AD_{min},GABA_{max})$. We choose small values of $AD_{min} = GABA_{min} = 0.001$ rather than exactly zero to prevent inclusion of a potentially unstable steady state. We can scale the size of this box according to the values of $AD_{max}$ and $GABA_{min}$. These two points are represented by the intersections of the $AD$ and $GABA_{VLPO}$ nullclines. To describe the activity between these two points, we begin with the fact that when $AD$ is at it's maximum value, $GABA_{VLPO}$ is at it's minimum value (0.001). This signifies that at the peak of wake activity, sleep activity is at its minimum. As we move along the $AD$ nullcline, the concentration of AD decreases to 0.001 as the concentration of $GABA_{VLPO}$ increases to $GABA_{max}$. From the point $(GABA_{VLPO}, AD) = (GABA_{max},0.001)$, we move along the $GABA_{VLPO}$ nullcline slowly until we reach the point $(GABA_{VLPO}, AD) = (0.001, AD_{max})$ again, creating a closed orbit. Figure \ref{fig:fig8} shows this movement in the phase plane, as well as showing the location of the fixed point and neighboring trajectories. Everything is in non-dimensional form, and therefore should be thought of in terms of threshold levels rather than micro-moles per liter. \\

The system is made up of two widely separated time scales: the slow, 24 hour time scale that is governed by the build-up and decay of adenosine, and the fast, O(1) (minute) time scale that involves interactions of all the sleep and wake active neutrotransmitters, but is regulated by the oscillation of the slow time scale. \\

We begin by forming the slow time system that will later be shown to be responsible for the limit cycle dynamics of the system. This system is comprised of two equations: the adenosine equation that fluctuates with wake activity, and the GABA in the VLPO equation that fluctuates with NREM sleep activity. We consider NREM sleep to be the ground state of the sleep-wake system, and all other oscillations to be deviations thereof. The adenosine equation must involve a constant source term, due to adenosine's constant creation as a by-product of neuronal activity, as well as a decay term that is dependent on adenosine's natural rate of degradation as well as it's loss from the synaptic gap through the its disinhibition of GABA in the VLPO \cite{adenosine_effects}. It is well known that mechanism for receptor binding is nonlinear \cite{neurotransport}, and so the disinhibition effect of $GABA_{VLPO}$ via adenosine should be nonlinear. It also follows that this nonlinear term should appear in both the equations for $AD$ and $GABA_{VLPO}$ because adenosine's disinhibition of the VLPO both depletes the concentration of adenosine while simultaneously (and with the same rate) adding to the $GABA_{VLPO}$ concentration. The $GABA_{VLPO}$ equation must contain the uptake term for GABA, a degradation term and a contribution from adenosine as well as the inhibitory and excitatory of other neurotransmitters. A constant loss term is required in the $GABA_{VLPO}$ equation as a representation of its quick breakdown time in the synaptic gap \cite{Winsky}. Though it is being broken down quickly and continuously, $GABA_{VLPO}$ still provides a tonic effect on the sleep-wake system as a whole \cite{Winsky}. Based on these requirements as well as the observations of the glycolysis reaction \cite{Selkov}, we obtain the following slow time system: \\

\noindent \textbf{Slow Time System}
\begin{equation}
\frac{dAD}{dt} = k_{1} - k_{2}AD - (GABA_{VLPO})^{2}(AD) 
\label{eq:eq11}
\end{equation}

\begin{equation}
\frac{dGABA_{VLPO}}{dt} = -\epsilon - k_{3}GABA_{VLPO} + k_{4}AD + (GABA_{VLPO})^{2}(AD) \label{eq:eq12}
\end{equation}

\bigskip\noindent
where
$$\epsilon = \mu - (a_{1}ACh_{BF} + a_{2}ACh_{LDT/PPT} + a_{3}NA + a_{4}S + a_{5}DA)$$ 
To provide a more physical explanation of the assumed reaction system, we know that adenosine is a bi-product of the constant source of energy in the nervous system, adenosine triphosphate (ATP). We know that adenosine is continuously converted back into ATP \cite{metabo}.  From experiments in the VLPO, we know that activation of VLPO neurons through adenosine-mediated inhibition contributes to the somnegenic actions of adenosine, and only adenosine acts directly on the VLPO in the sleep-wake system \cite{adenosine_effects,Arrigoni,adenosine1,adenosine2}. During wake, increased levels of adenosine promote sleep by disinhibiting the release of GABA in VLPO. Adenosine levels in this area increase during wake and decrease during sleep \cite{adenosine_effects}. In this way the physical foundation for the glycolysis system is similar to and provides a good framework for our slow time scale system. \\

\subsubsection{Characterizing REM/NREM}

REM and NREM sleep are determined by a nonlinear combination of neurotransmitters.  Specifically, $GABA_{VLPO}$ promotes the state of REM sleep, as does acetylcholine from both the LDT/PPT and the BF. While NREM can be seen as an absence of the vigilance state and the contributions from GABA, REM sleep must incorporate contributions from neurotransmitters active in the wake state as well as those in the NREM state. Letting $R(t)$ represent the REM sleep state, we assume that $R'' + \Phi R'+ \Gamma R= 0$. To determine  $\Phi$, note that the transition from the waking state to the sleeping state, and vice versa, occurs when $GABA_{VLPO} = AD$.  Given that the REM state has an asymptotically stable steady state during wake, which occurs when  $AD>GABA_{VLPO}$, and has a limit cycle during sleep, which occurs when $AD<GABA_{VLPO}$, then $\Phi$  contains a term of the form $(GABA_{VLPO} - AD)(R^2-R_0^2)$.  When entering wake, so $GABA_{VLPO}=AD$ and $\dot{AD}>0$, there is an impetus for the impending dormancy of REM and the resulting mathematical requirement is that $\Phi(\dot{AD},AD)>0$. Likewise, when entering sleep, so $GABA_{VLPO}=AD$ and $\dot{AD}<0$, the requirement for the oncoming REM limit cycle is $\Phi(\dot{AD},AD)<0$. This gives us

\begin{equation}
\frac{d^2 R}{dt^2} + \left [\alpha (GABA_{VLPO} - AD)(R^2 - R_0^2) + \beta \frac{dAD}{dt} \right ] \! \frac{dR}{dt} + \Gamma R = 0 \,. \label{eq:eq14}
\end{equation}

\subsubsection{Some Special Cases}
Given the complexity of (\ref{eq:eq14}), it is worth considering the solution in particular situations and then using this to help estimate the oscillatory response of equation (\ref{eq:eq14}).  For example,  the period for $R$ can be determined from the fact that there  are  three to four normal REM oscillations  over a normal night's sleep.   So, when $GABA_{VLPO} = AD$, 

\begin{equation}
\frac{d^2 R}{dt^2} + \beta \frac{dAD}{dt}\frac{dR}{dt} + \Gamma R = 0 \,. \label{eq:eq15}
\end{equation}

\noindent 
Given the slowly varying nature of $AD$, then $R = e^{rt}$, where $r =\frac{1}{2}( -\beta u \pm \sqrt{\beta^2 u^2 - 4\Gamma})$ and $u=(AD)'$.  Since $u$ strongly controls the decay, and to have $\Gamma$ control the number of REM oscillations, we require $r = -u \pm i \sqrt{\gamma}$, where $\sqrt{\gamma}$ is approximately the number of oscillations. Thus, $\beta = 2$ and $\Gamma = \gamma + u^{2}$.  From this we have that 
\begin{equation}
\frac{d^2 R}{dt^2} + 2 \! \left [(GABA_{VLPO} - AD)(R^2 - 1.3) + \alpha \frac{dAD}{dt} \right ] \!  \frac{dR}{dt} +   \! \left (\gamma + \left (\frac{dAD}{dt}\right)^2 \right )  \! R= 0 \, . \label{eq:eq16}
\end{equation}

\noindent 
where $\alpha $ is determined numerically and $8< \gamma<16$. \\

Based on Lienard's Theorem, we know that the above oscillator has a unique, stable limit cycle \cite{strogatz1}. Thus the REM oscillations will continue to function in a stable, periodic fashion even when subjected to reasonable perturbations. Evidence of this result is in the Orexin Knockout section. \\

\subsection{Analysis}
To find the values of the coefficients of the system and to ensure stability, we present several methods of analyzing the  equations. We begin with the limit cycle coming from the Slow Time System.   As stated earlier, the adenosine equation drives the wake state dynamics, and the GABA equation drives the sleep state dynamics.  The specific equations are

\begin{equation}
\frac{d}{dt} AD = k_{1} - k_{2}AD - (GABA_{VLPO})^{2}(AD) \label{eq:eq17}
\end{equation}
\begin{equation}
\frac{d}{dt} GABA_{VLPO} = -\epsilon - k_{3}GABA_{VLPO} + k_{4}AD + (GABA_{VLPO})^{2}(AD) ,
\label{eq:eq18}
\end{equation}

\noindent
where
$$\epsilon = \mu - (a_{1}ACh_{BF} + a_{2}ACh_{LDT/PPT} + a_{3}NA + a_{4}S + a_{5}DA) .$$

It is not possible to determine the rate constants with the available experimental evidence.  So, we take a different tack and establish bounds on the constants that will produce the responses known to occur in sleep and wake.  To begin, we compute the nullclines of the slow time system in order to create a trapping region for the limit cycle and to determine the coefficients of the system by using an approximate fixed point. We have four coefficients to determine ($k_1$, $k_2$, $k_3$, $k_4$) as well as reasonable values for the parameter $\epsilon$. Plugging the approximate fixed point into both nullclines will give us two equations with five unknowns. In order to make the system well determined, we observe the maxima and minima of the equations as well as the linearized slow time system with the approximate fixed point. Knowing the form of the determinant and the trace of this system needed to maintain stability, we can determine the rest of the unknowns. \\

Now, solving the nullcline equations obtained from \ref{eq:eq17} and (\ref{eq:eq18}) for $AD$ and we have
$$AD = \frac{k_{1}}{k_{2}+GABA_{VLPO}^2} \, ,$$
and 
$$AD = \frac{k_{3}GABA_{VLPO} - \epsilon}{k_{4}+(GABA_{VLPO})^2} \, .$$

\noindent
Setting these equal to each other, we have three points of intersection, only one of which is real.  Using $GABA_{VLPO} \in [0,2]$ and $AD \in [0,2]$ as the boundaries of  the trapping region, we can approximate a fixed point inside the closed orbit at $GABA_{VLPO} = 0.703$ and $AD = 0.823$.  This provides us with two of the needed equations.\\

We also know that each of these equations is zero separately when $AD = AD_{max}$ and $GABA_{VLPO} = GABA_{min}$, as well as when $AD = AD_{min}$ and $GABA_{VLPO} = GABA_{max}$. We let $AD_{max} = GABA_{max} = 2$ and $AD_{min} = GABA_{min} = 0.01$.  Taking $\epsilon = 0.3$, and using the stated values for the max and min concentrations it follows that 
$k_{1} = 0.49$, 
$k_{2} = 0.1$, 
$k_{3} = 0.3$. 
$k_{4} = 0.15$.
 We also note that $0.29 < \epsilon \leq 0.32 $ for the system to remain stable. \\

\subsubsection{Remaining Coefficient Estimation}

Now that the coefficients for the slow time system are determined, we can use the requirements that the eigenvalues have either pure imaginary form or negative real parts in the complex form to determine the coefficients of the entire system. We have the solution vector of the linearized system as $\mathbf{x} = 
[ GABA_{BFw},
GABA_{BFs},
OX,
H,
ACh_{BF} ,$ 
$ACh_{LDT/PPT} ,
NA ,
S ,
DA ,
AD ,
GABA_{VLPO} 
]^T$. The matrix of coefficients and linearized values for the system are  given in Table \ref{tb:1.euler}.

\begin{sidewaystable}
\centering
\vspace{1in}
\renewcommand\arraystretch{1.5}
\renewcommand\tabcolsep{2pt}
B =
\begin{tabular}{p{0.5in} p{0.5in} p{0.5in} p{0.5in} p{0.5in} p{0.5in} p{0.7in} p{0.7in} p{0.6in} p{0.6in} p{1in}}
\footnotesize{-$c_{1}$-GA1} & 0 & \footnotesize{$c_{2}$} & 0 & \footnotesize{$c_{3}$} & 0 & 0 & 0 & 0 & \footnotesize{$c_{4}$} & 0 \\
0 & \footnotesize{-$c_{5}$-GA1} & 0 & 0 & 0 & 0 & -$c_{6}$ & 0 & 0 & -$c_{7}$ & $c_{8}$ \\
0 & -$c_{9}$ & \footnotesize{$c_{10}$-$c_{11}$} & 0 & 0 & $c_{12}$ & \footnotesize{$c_{13}$-$c_{14}$} & -$c_{15}$ & 0 &  $c_{16}$ & $-c_{17}$ \\
0 & 0 & $c_{18}$ & \footnotesize{$c_{19}$-$c_{20}$} & $c_{21}$ & 0 & $c_{22}$ & 0 & 0 & $c_{23}$ & -$c_{24}$ \\
0 & -$c_{25}$ & $c_{26}$ & $c_{27}$ & -$c_{28}$ & 0 & 0 & -$c_{29}$ & 0 & $c_{30}$ & 0 \\
0 & 0 & $c_{31}$ & $c_{32}$ & 0 & -$c_{33}$ & -$c_{34}$ & -$c_{35}$ & 0 & $c_{36}$ & 0 \\
0 & 0 & $c_{37}$ & 0 & $c_{38}$ & 0 & \footnotesize{-$c_{39}$-$c_{40}$-hN} & 0 & 0 & $c_{41}$ & -$c_{42}$ \\
0 & 0 & $c_{43}$& $c_{44}$ & 0 & 0 & $c_{45}$ & \footnotesize{-$c_{46}$-$c_{47}$-hS} & 0 & $c_{48}$ & -$c_{49}$ \\
0 & 0 & $c_{50}$ & 0 & 0 & $c_{51}$ & -$c_{52}$ & $c_{53}$ & \footnotesize{-$c_{54}$-$hD$} & $c_{55}$ & -$c_{56}$ \\
0 & 0 & 0 & 0 & 0 & 0 & 0 & 0 & 0 & \footnotesize{-$0.1 - G_{V}^2$} & \footnotesize{-$2G_{V}A$}\\
0 & 0 & 0 & 0 & \footnotesize{-$c_{57}$} & \footnotesize{-$c_{58}$} & \footnotesize{-$c_{59}$} & \footnotesize{-$c_{60}$} & \footnotesize{-$c_{61}$} &  \footnotesize{ $0.15 + G_{V}^2$} & \footnotesize{-$0.3+2G_{V}A$}
\end{tabular}
\caption{Matrix B of the linearized system. Note that  $hN=hNET$, $hS=hSERT$, $hD=hDAT$, $G_{V} = GABA_{VLPO}$ and $A = AD$.}
\label{tb:1.euler}
\end{sidewaystable}

In order for the system to be stable, we require that the real parts of the eigenvalues to be negative. Restricting the coefficients to all be between 0 and 1, we use an algorithm that generates coefficients that keep the real parts of the eigenvalues negative. We also note that the coefficients in $GABA_{VLPO}$ equation need to be small compared to the rest of the coefficients in order to keep $\epsilon$ bounded.

\begin{table}
\begin{tabular}{p{1in} p{1.2in} p{1in} p{1in}}
\hline 
Parameter & Value & Parameter & Value\\
\hline
$GA1$ & 1 & $c_{27}$ & 0.10973\\
$hNET$ & 0.457 & $c_{28}$ & 0.32943\\
$hSERT$ & 0.463 & $c_{29}$ & 0.57879\\
$hDAT$ & 1.22 & $c_{30}$ & 1\\
$k_{1}$ & 0.49 & $c_{31}$ & 0.02091\\
$k_{2}$ & 0.1 & $c_{32}$ & 0.12648\\
$k_{3}$ & 0.3 & $c_{33}$ & 0.23472\\
$k_{4}$ & 0.15 & $c_{34}$ & 0.57122\\
$\mu$ & 0.3 & $c_{35}$ & 0.02332\\
$\alpha$ & 1 & $c_{36}$ & 1\\
$\gamma$ & 8 & $c_{37}$ & 0.61305\\
$c_{1}$ & 0.75709 & $c_{38}$ & 0.06864\\
$c_{2}$ & 0.28014 & $c_{39}$ & 0.08638\\
$c_{3}$ & 0.61048 & $c_{40}$ & 0.1\\
$c_{4}$ & 0.76636 & $c_{41}$ & 0.543\\
$c_{5}$ & 0.32431 & $c_{42}$ & 0\\
$c_{6}$ & 0.83153 & $c_{43}$ & 0.13822\\
$c_{7}$ & 0.03471 & $c_{44}$ & 0.35956\\
$c_{8}$ & 0.01 & $c_{45}$ & 0.11839\\
$c_{9}$ & 0.19577 & $c_{46}$ & 0.13753\\
$c_{10}$ & 0.79157 & $c_{47}$ &0.1\\
$c_{11}$ & 0.97026 &  $c_{48}$ & 0.537\\
$c_{12}$ &-1 &  $c_{49}$ & 0\\
$c_{13}$ & 0 & $c_{50}$ & 0.67749\\
$c_{14}$ & 0.36341 & $c_{51}$ & 0.53609\\
$c_{15}$ & 0.70633 & $c_{52}$ & 0.36464\\
$c_{16}$ & 0.1 & $c_{53}$ & 0.59591\\
$c_{17}$ & 0.97643 & $c_{54}$ & 0.75091\\
$c_{18}$ & 0.56740 & $c_{55}$ & 0.1\\
$c_{19}$ & 0.91859 & $c_{56}$ & 0.22\\
$c_{20}$ & 0.50364 & $c_{57}$ &$10^{-5}$\\
$c_{21}$ & 0.23758 & $c_{58}$ & $10^{-5}$\\
$c_{22}$ & 0.1 & $c_{59}$ & $10^{-5}$\\
$c_{23}$ & 1 & $c_{60}$ &$10^{-5}$\\
$c_{24}$ & 0 & $c_{61}$ &$10^{-5}$ \\
$c_{25}$ & 0.92037 & & \\
$c_{26}$ & 0.04185 & & \\

\hline
\end{tabular}
\caption{Coefficient and parameter values used in the sleep-wake cycle model.}
\end{table}

\subsubsection{Sub-Systems}
When analyzing such a large system of equations with undetermined coefficients, it is helpful to break the large system into smaller sub-systems and observe the dynamics on a local scale. We can utilize the neurophysiology of the physical neurotransmitter sub-systems to determine the coefficients as well as the stability profiles of the individual sub-systems. We break the system up into to following categories:
NREM Specific,
Sleep Specific (NREM Inclusive),
Wake Specific,
REM + NREM Specific,
REM + Wake Specific,
Wake Inclusive,
Sleep Inclusive,
REM Inclusive.
Note that there is no set of equations associated with REM only. The categories represent what neuron populations are active during the corresponding states. We observe the steady-state behavior of each sub-system. Since we know the relative magnitude of the neurotransmitter amounts as well as their dynamic behavior during the states of REM, NREM and Wake, the coefficients of these smaller, more manageable systems should be more easily attainable. This also allows us to observe the more manageable systems for potentially extreme behavior, and provides us with bounds to prevent this from happening. This analysis does not give us specific values for the system coefficients, knowing the relative values makes visualizing the system easier as well as simplifying the eigenvalue analysis. \\

\subsubsection*{NREM Specific}
The NREM specific system involves the neuron population only active during the state of NREM sleep. The neuron population that is only active in NREM sleep is the GABAergic population in the BF. We assume all equations but $\frac{dGABA_{BFs}}{dt}$ to be equal to zero, making each neurotransmitter other than $GABA_{BFs}$ constant, and assume the other neurotransmitters in the equation to be at their steady state values when only $GABA_{BFs}$ is active. This is not necessarily a physiologically accurate view of the system for an extended length of time, but we can view a snapshot of this activity. We follow this same criteria for the rest of the small systems.

\subsubsection*{Wake Specific}
The Wake specific system involves neuron populations with the most activity during the wake state. The monoaminergic neuron populations (including OX, H, NA, S and DA) are most active in this state.

\subsubsection*{Sleep Specific (NREM Inclusive)}
The Sleep Specific system involves neurotransmitters that promote sleep states (REMS and NREMS). The populations that are active in this state are GABAergic populations in the VLPO and BF. 

\subsubsection*{REM + Wake Specific}
The REM and Wake Specific system involves neurotransmitters that promote REMS and the Wake state. The neuron populations that are active during these states are GABAergic populations in the BF and cholinergic populations in the BF and LDT/PPT.

\subsubsection*{REM + NREM Specific}
The REM and NREM Specific system involves neurotransmitters promoting REM sleep and the NREM sleep. The neuron population that is active during these states is the GABAergic population in the VLPO. 

\subsubsection*{Wake Inclusive}
The Wake Inclusive system involves neurotransmitters that promote the Wake state as well as those that promote both Wake and REMS. The neuron populations that are active during these states are the GABAergic population in the BF, the cholinergic population in the BF and LDT/PPT, and the monoaminergic populations. 

\subsubsection*{Sleep Inclusive}
The Sleep inclusive system involves neurotransmitters associated with NREM specific, REM and NREM specific, and Wake and REM specific states. The neuron populations that are active during these states are GABAergic populations in the VLPO and BF, and cholinergic populations in the BF and LDT/PPT.

\subsubsection*{REM Inclusive}
The REM Inclusive system involves neurotransmitters associated with the Wake and REM specific and REM and NREM specific states. The neuron populations that are active during these states are the GABAergic population in the VLPO, the wake promoting GABAergic population in the BF as well as the cholinergic populations in the BF and LDT/PPT. \\

This analysis allows for an approximation of the coefficients, but it does not yield exact values for these coefficients. Figures \ref{fig:fig12} and \ref{fig:fig13} show the change in concentration levels of Acetylcholine in the BF and GABA in the BF, respectively over a period of 24 hours. The remainder of the neurotransmitters involved show the same fluctuations as that of Acetylcholine, whereas GABA in the sleep active portion of the BF is unique. \\

\begin{SCfigure}[][h]
\centering
\includegraphics[width=2in]{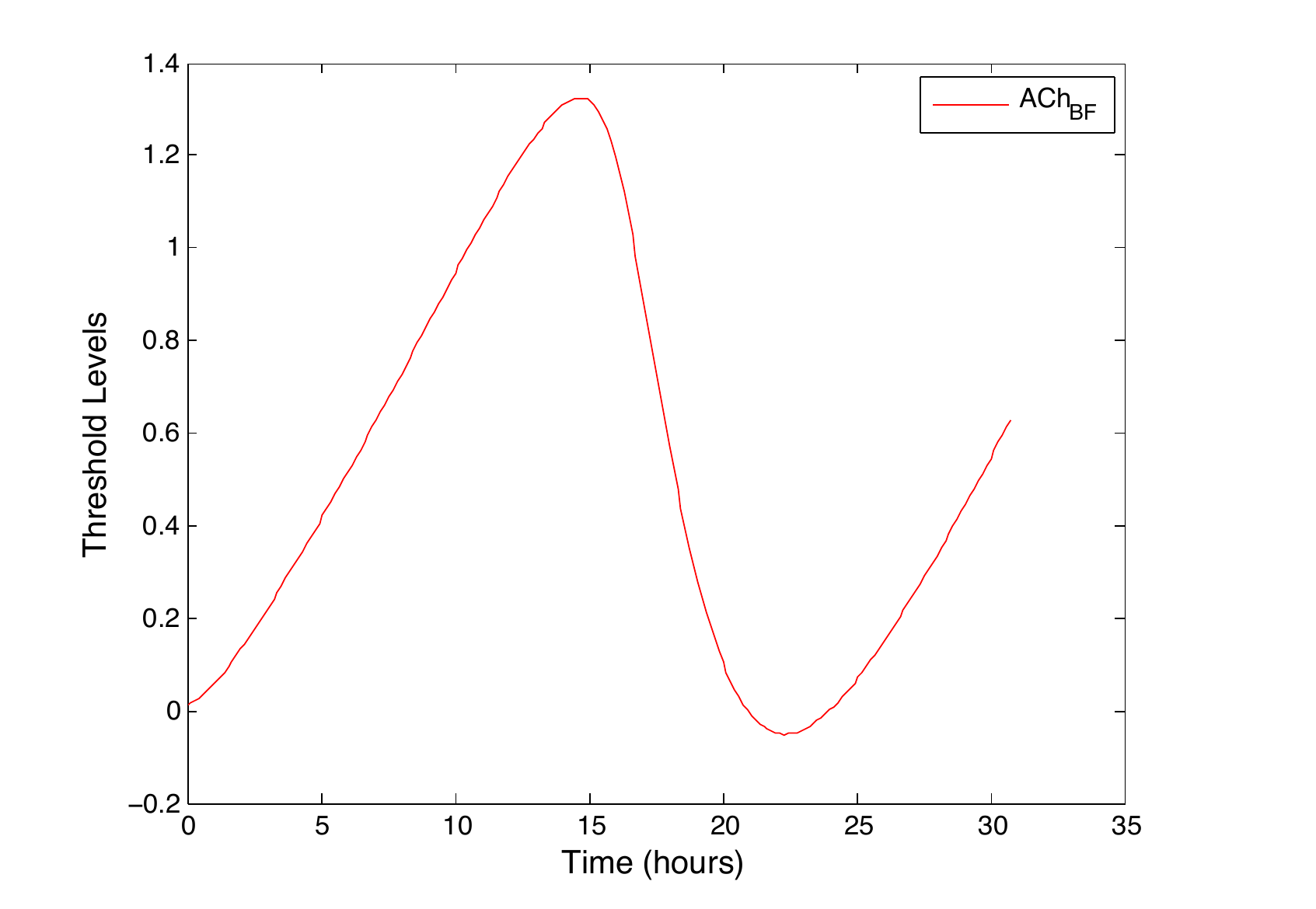}
\caption{Concentration levels of Acetylcholine in the BF over 35 hours.}
\label{fig:fig10}
\end{SCfigure}

\begin{SCfigure}
	\centering
		\includegraphics[width=2in]{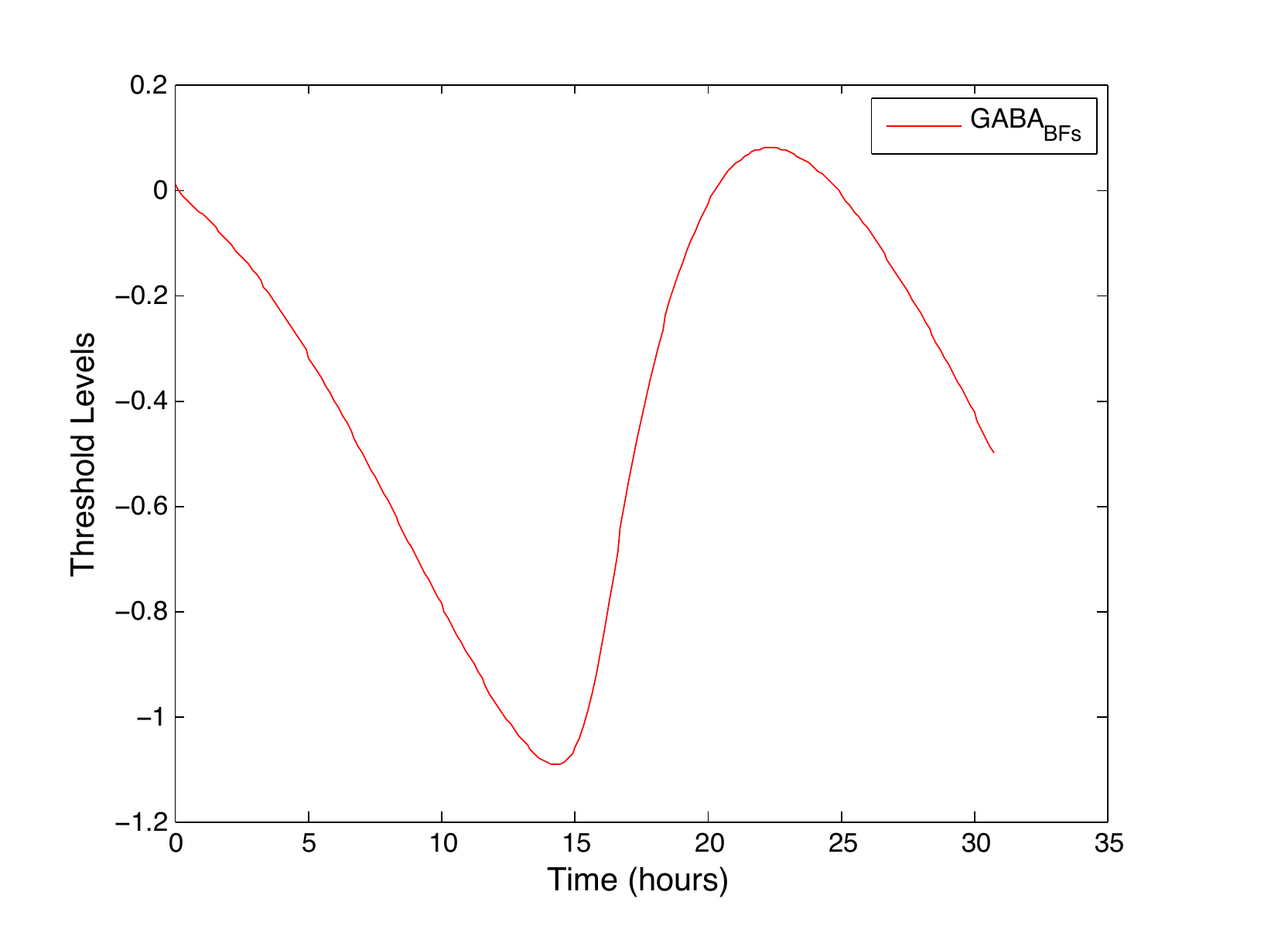}
	\caption{Concentration levels of GABA in the BF associated with sleep over 35 hours.}
	\label{fig:fig11}
\end{SCfigure}

\section{Results and Discussion}

\begin{figure}[t]
	\begin{center}
		\includegraphics[width=4.5in]{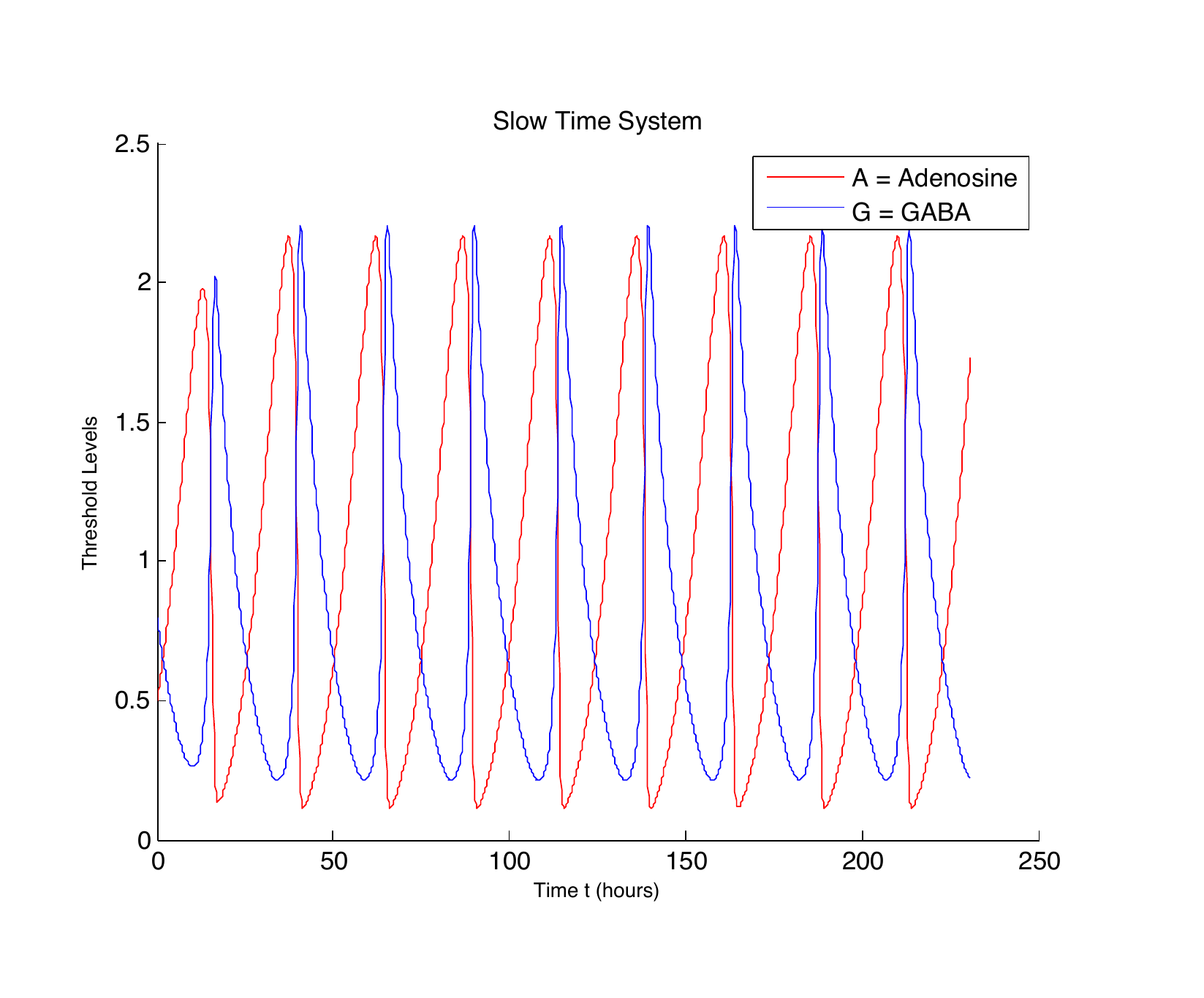}
	\caption{Adenosine and GABA (in the VLPO) levels of the sleep-wake cycle model over 216 hours (9 days).}
	\label{fig:fig12}
 \end{center}
\end{figure}

\begin{figure}[t]
	\begin{center}
		\includegraphics[width=5in]{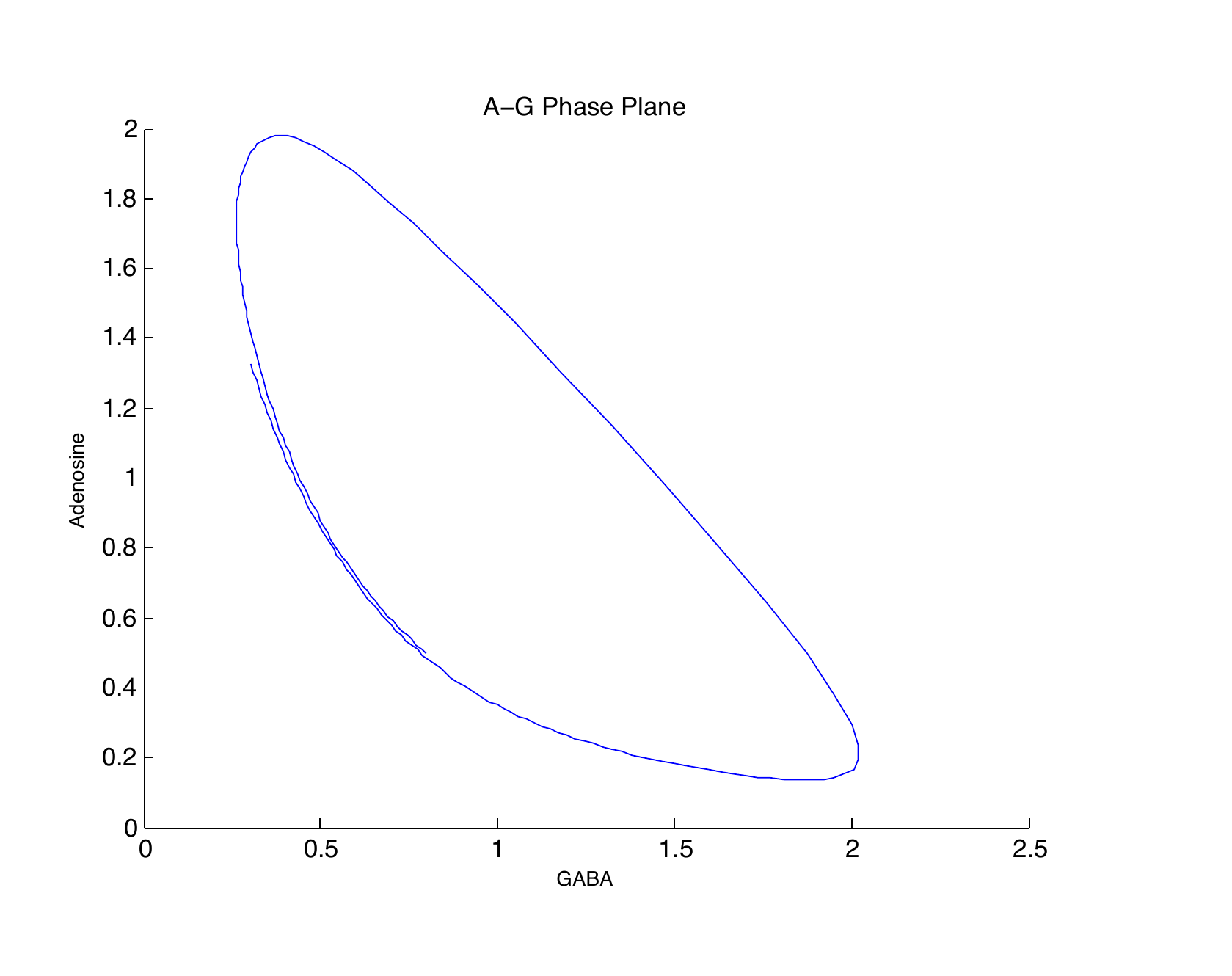}
	\caption{Phase plane of sleep-wake cycle model for one cycle.}
	\label{fig:fig13}
 \end{center}
\end{figure}

\begin{figure}[t]
	\begin{center}
		\includegraphics[width=5in]{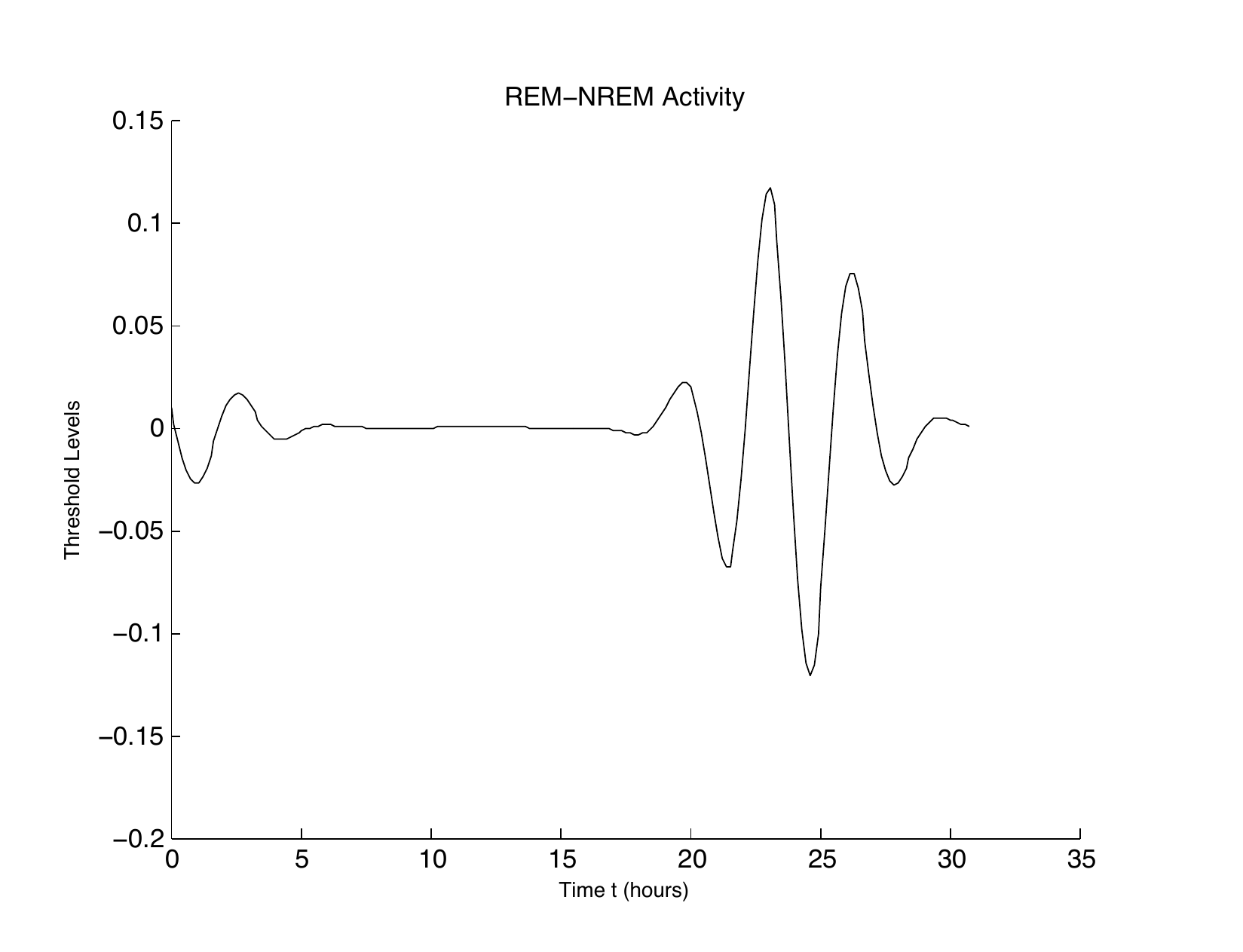}
	\caption{REM cycle visualization of sleep-wake cycle model over a day.}
	\label{fig:fig14}
 \end{center}
\end{figure}

\begin{figure}[t]
	\begin{center}
		\includegraphics[width=5in]{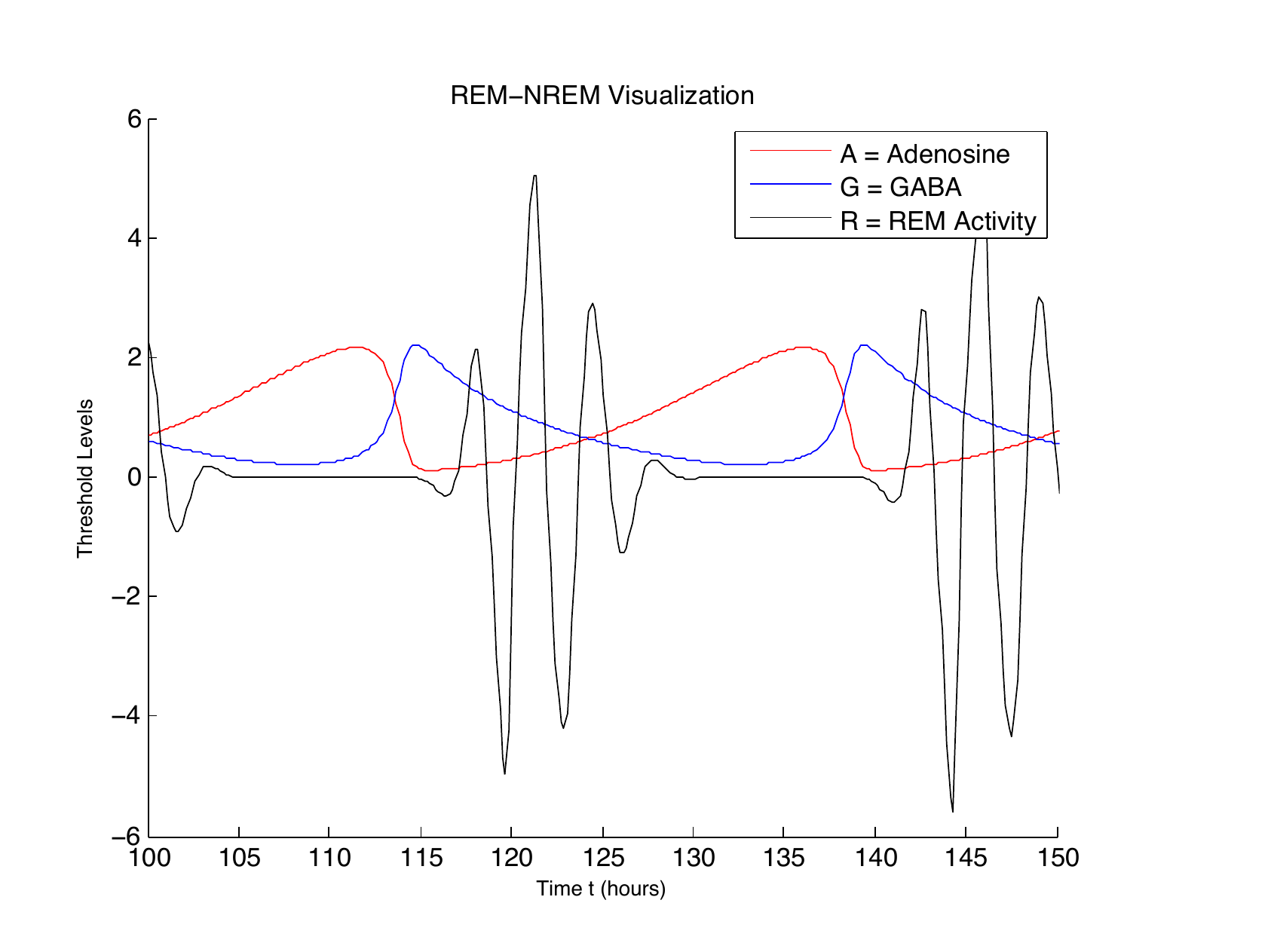}
	\caption{Unperturbed system of AD, GABA and REM from hours 150 to 200 - over two days. }
	\label{fig:fig15}
 \end{center}
\end{figure}

The results of the implemented normal, unperturbed, system are shown in Figures \ref{fig:fig12}-\ref{fig:fig15}. Figure \ref{fig:fig12} depicts the slow time variables (which are the easiest to visualize on a large time scale) $AD$ and $GABA_{VLPO}$. The $AD$ curve represents the non-dimensional amount of adenosine accumulation in the BF, which consequently measures sleep pressure. The $GABA$ curve represents the non-dimensional amount of GABA accumulation in the VLPO, which is a measure of wake pressure. Both exhibit a rapid transition between sleep and wake states as well as ongoing periodic behavior for nine days. The system runs like this for much larger periods, demostrating that unless perturbed, the sleep-wake cycle is a self-sustained, periodic oscillation. The period of each cycle is approximately 24 hours - the wake state occupies roughly two-thirds of this time, while the sleep state makes up the remaining one-third. This is consistent with our knowledge of the structure of the human sleep-wake cycle. Figure \ref{fig:fig13} depicts one cycle in the $GABA-AD$ phase plane, showing the exact behavior we required from our previous nullcline analysis. Figure \ref{fig:fig14} is the isolated REM curve for 30 hours. Figure \ref{fig:fig15} shows the $AD$, $GABA$ and REM curves simultaneously in a sample 48 hour time slice. The REM/NREM oscillations respond quickly to the transitions between wake and sleep, and cycle through three to four oscillations in the designated eight hour sleep period. From the biological background and neurotransmitter details, this displays physiologically reasonable results.\\

\subsection{Orexin Knockout Experiment}

Narcolepsy is a neurodegenerative disease where there is a loss of neurons containing orexin (OX) \cite{CSF}. It has been shown to exist across many mammals, and has been studied extensively in mice. Here, we take information from the sleep-wake cycles of orexin knock-out (which means the targeted destruction of orexin neurons) mice. These mice exhibit the narcolepsy symptoms of sleep fragmentation, sleep onset during REM sleep states and increased spontaneous REM and NREM sleep during the vigilance state \cite{OXKO,Wu,Baumann,Lu2,Grivel}. \\

\begin{figure}[h]
	\begin{center}
		\includegraphics[width=5in]{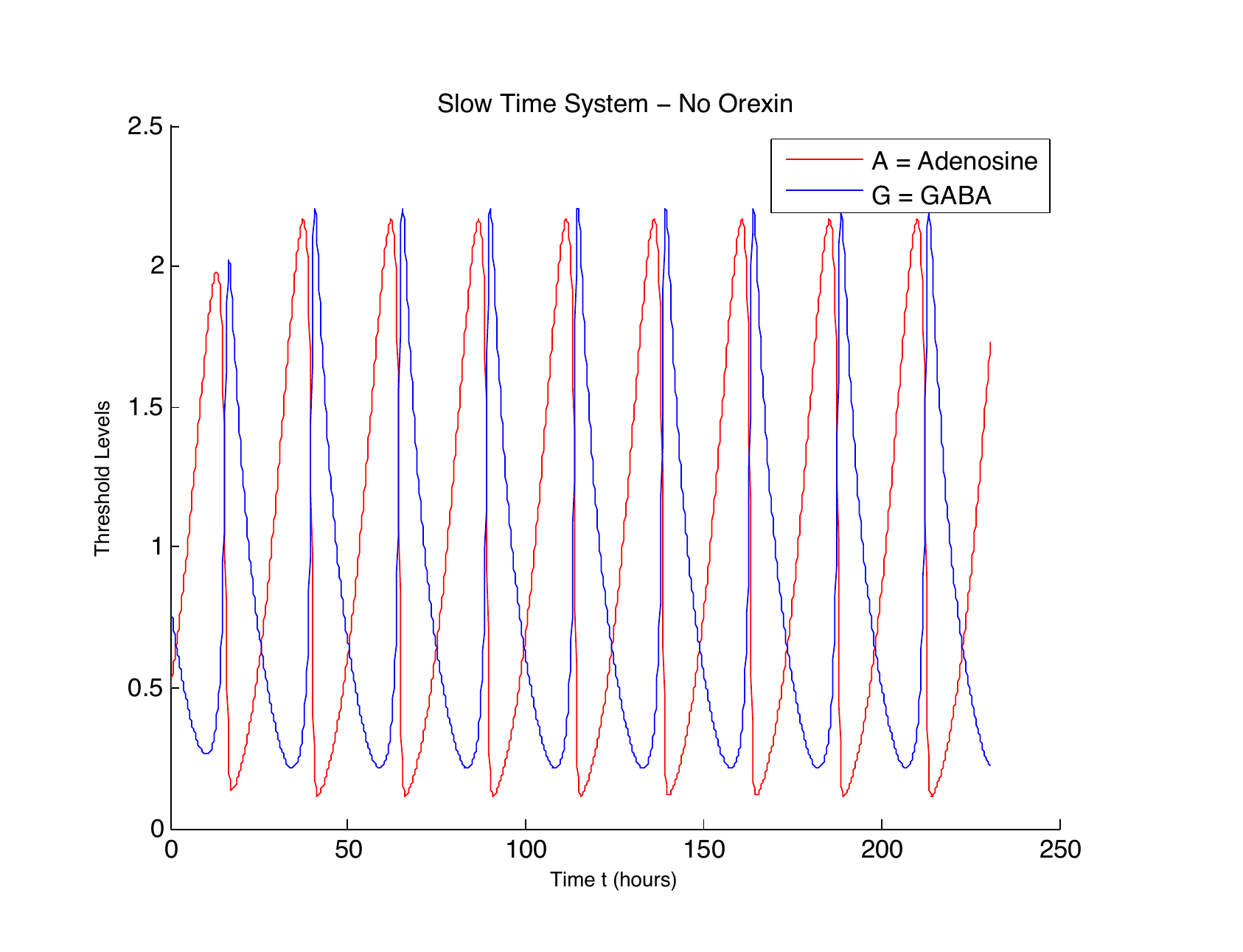}
	\caption{Orexin knockout system over 216 hours (9 days).}
	\label{fig:fig16}
 \end{center}
\end{figure}

A reasonable method to represent the orexin knock-out perturbation to the system is to change the flux of the OX amount equation. Rather than take orexin out all together, we limit the rate at which orexin is produced, but keep the removal factors the same. This is a reasonable way to represent the physiological destruction of orexin neurons, since none of the uptake and removal mechanisms are necessarily effected. It has been shown that levels of orexin decrease in human cerebral spinal fluid by at least 80 percent \cite{CSF}, and so we multiply orexin producing elements of the OX equation ($ACh_{LDT/PPT}$, $NA$, $OX$ and $AD$) by 0.2. Implementing these changes numerically, the periodicity and limit cycle of the system is maintained with orexin decrease. This is illustrated in Figure \ref{fig:fig16}, which clearly shows the same characteristics as seen in Figure  \ref{fig:fig12}.  However, we also don't see the stereotypical changes to the transition rates between behavioral states characteristic to narcolepsy. We have shown that the system can withstand perturbations \cite{OXKO}, but may need to think of implementing another method of making genetically specific perturbations. This also shows that orexin may have to play a larger role in measuring behavioral state instability, and may have to withstand more non-linear changes to its equation when investigated further. 

\subsection{Sleep Camp Data}
\cite{Carskadon} describes experiments related to perturbations in an individuals sleep/wake cycle.  The trials performed were carried out in a controlled situation with a group of 40 adolescents over a period of one week. The perturbations performed were delayed sleep and wake onset, prolonged sleep and wake periods and forced wake. These perturbations were achieved without the use of external chemical stimuli, and so were essentially achieved in a noninvasive fashion. Salivary cortisol and melatonin levels were recorded daily as well as frequent EEG readings for accurate brainwave activity. The initial time and length of the perturbation period were recorded, and so this data was easy to translate into perturbations in the model. About 50 perturbations were introduced (which translates to a period of about 10 days) and their effect on the system was recorded using graphical and data analysis methods. \\

\begin{figure}[!t]
	\begin{center}
		\includegraphics[width=5in]{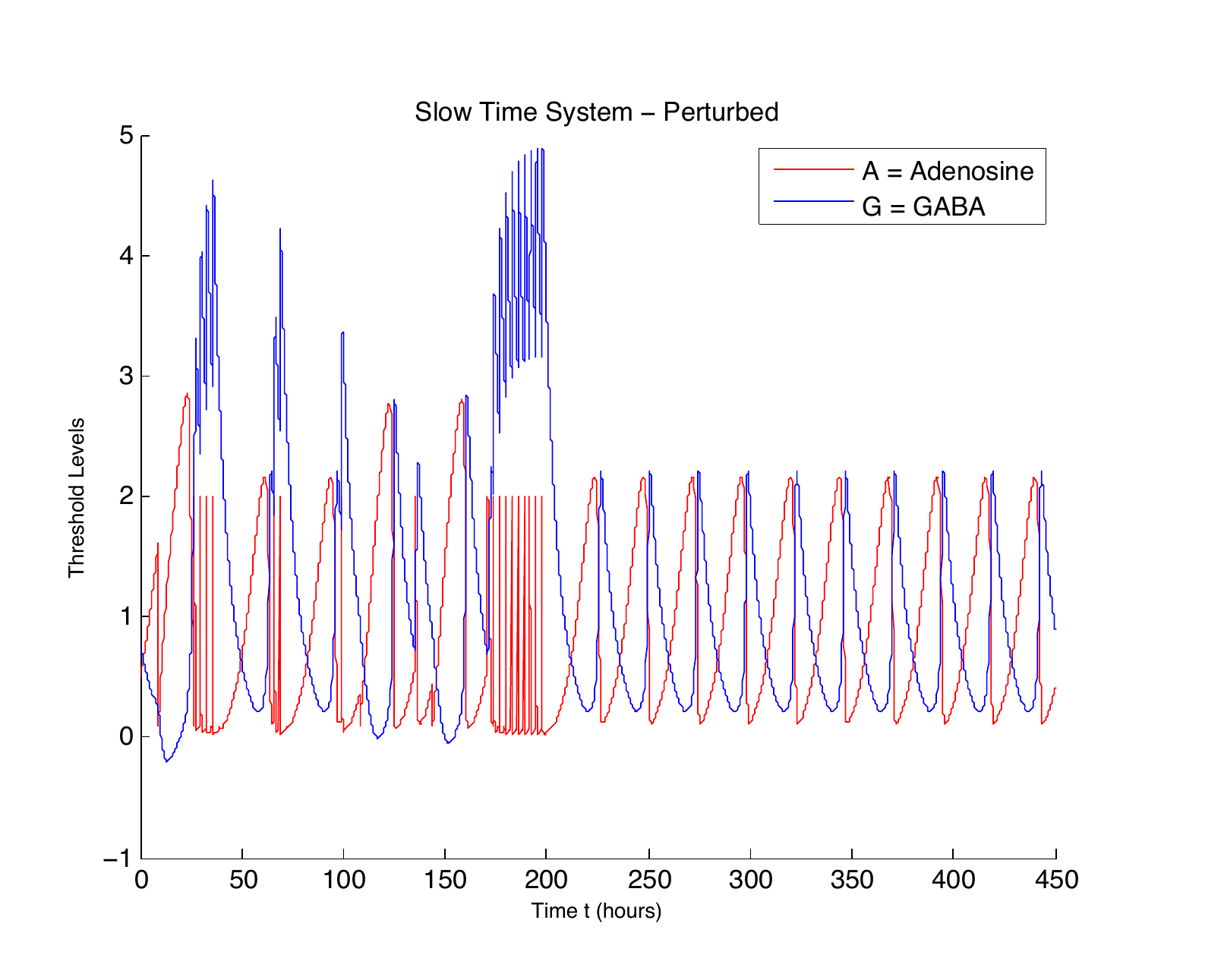}
	\caption{Perturbed slow time system with forced wake onset and forced sleep onset over a one week period. Extended time to show phase drift and maintained periodicity.}
	\label{fig:fig17}
 \end{center}
\end{figure}

\begin{figure}[t]
	\begin{center}
		\includegraphics[width=5in]{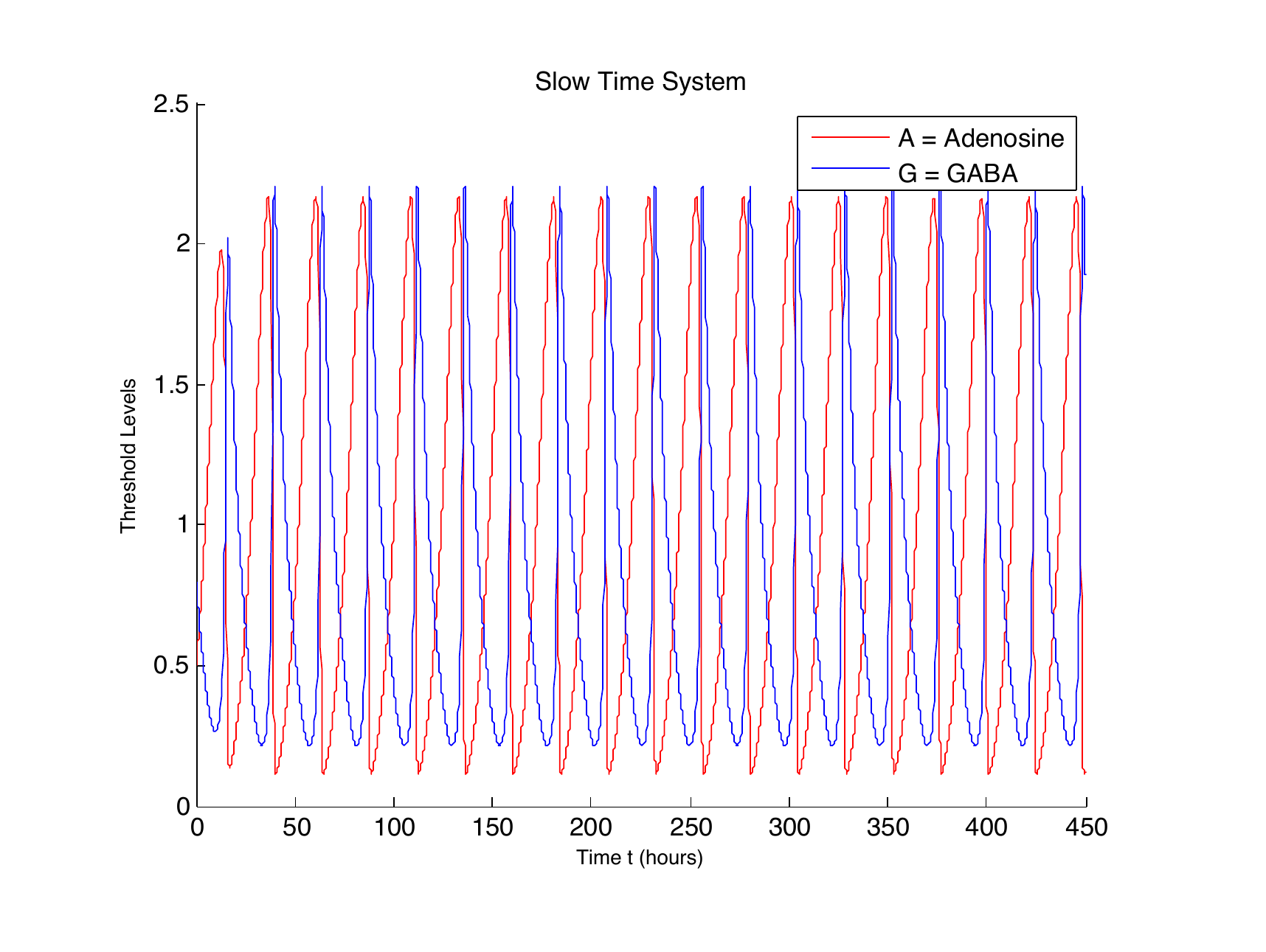}
	\caption{Unperturbed cycle allowed to run for the same amount of time as the perturbed cycle. }
	\label{fig:fig18}
 \end{center}
\end{figure}

\begin{figure}[h!]
	\begin{center}
		\includegraphics[width=5in]{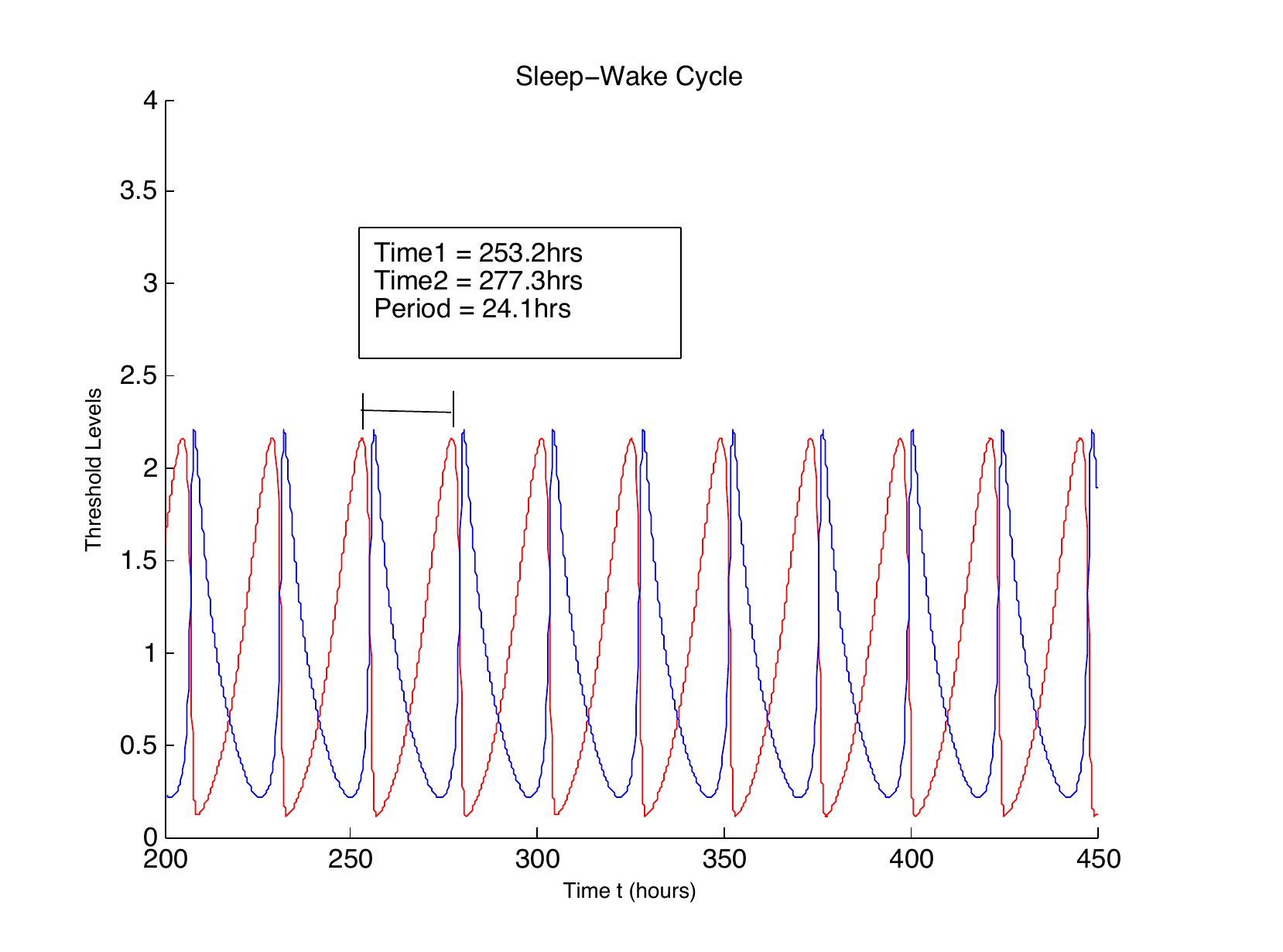}
	\caption{Examination of the length one period in an unperturbed cycle.}
	\label{fig:fig19}
 \end{center}
\end{figure}

\begin{figure}[t]
	\begin{center}
		\includegraphics[width=5in]{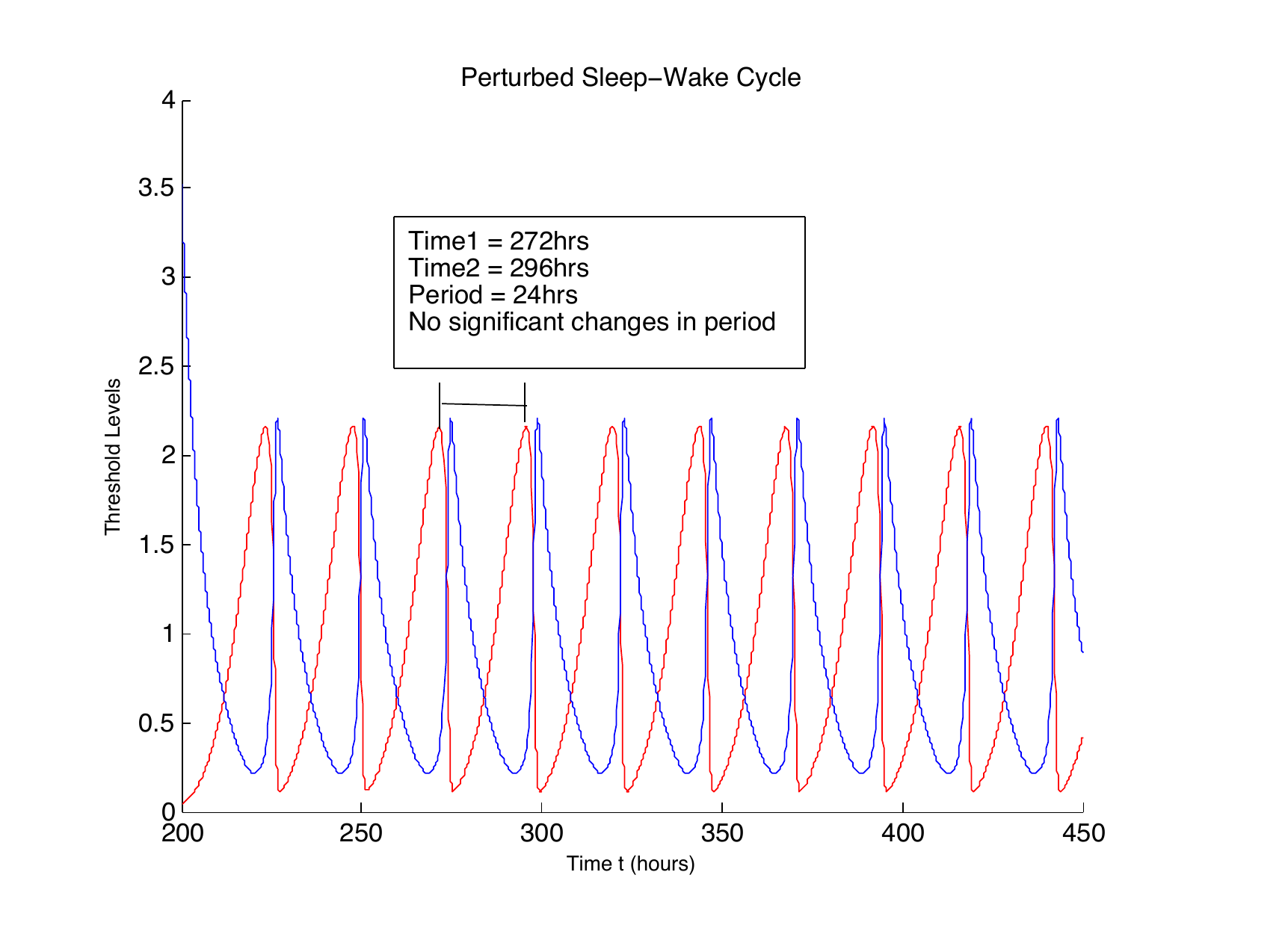}
	\caption{Examination of the length of one period in a perturbed cycle. }
	\label{fig:fig20}
 \end{center}
 
\end{figure}

\begin{figure}[b]
	\begin{center}
		\includegraphics[width=5in]{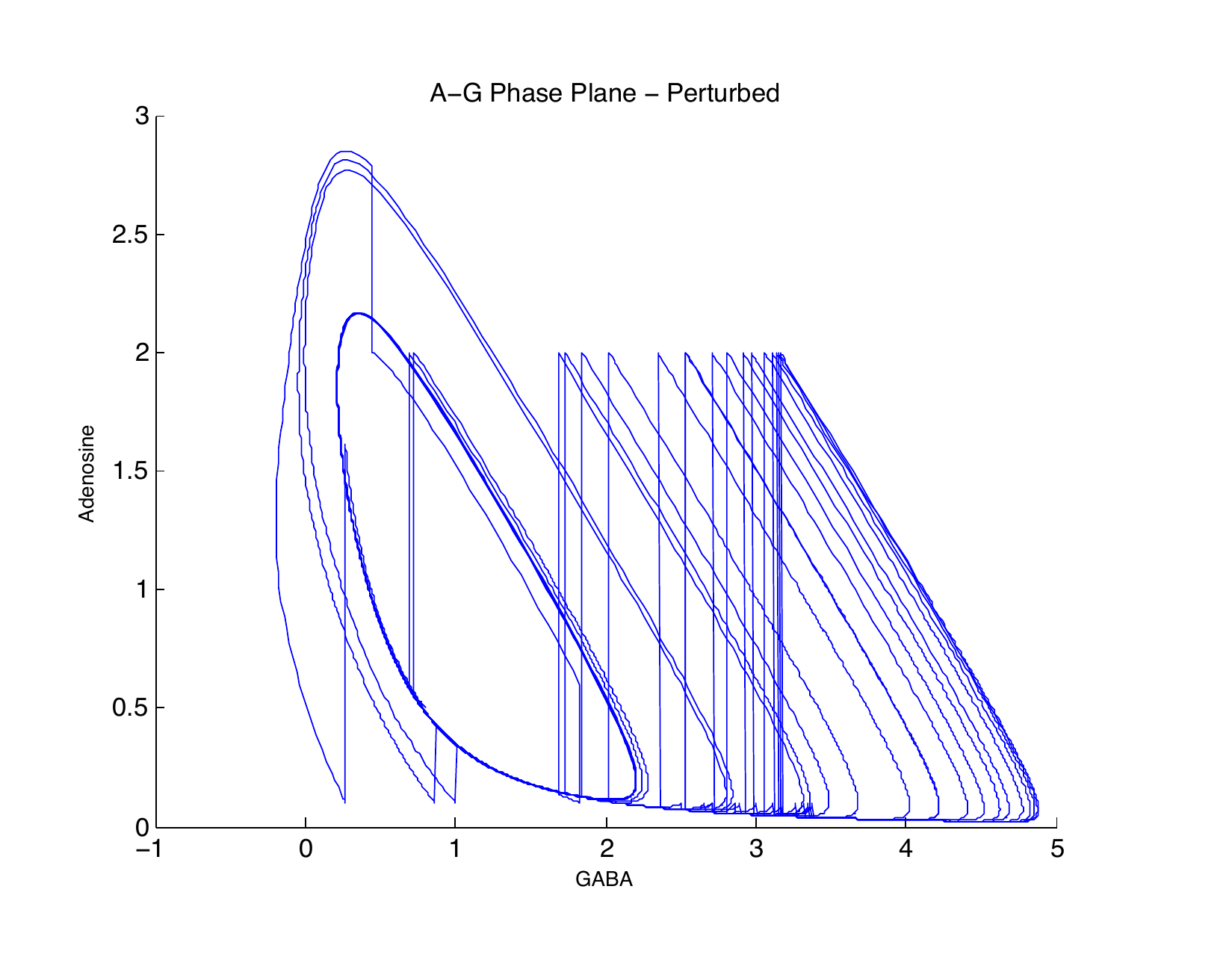}
	\caption{Perturbed GABA-Adenosine phase plane. }
	\label{fig:fig21}
 \end{center}
\end{figure}

The perturbations to the system were only introduced to one variable - adenosine. Forced sleep and wake were achieved by forcing $AD$ to its maximum at the appropriate forced waking initialization, and forcing $AD$ to its minimum at the appropriate forced sleep time. These sleep and wake times as well as the prolonged sleep and wake periods were taken from a data set of seven adolescent individuals. Once perturbed the amount of times done so in the physical experiment, we allow the system to run without perturbation for several cycles. This entire set of perturbations and the response is depicted in Figure \ref{fig:fig17}. Figure \ref{fig:fig18} shows what the subject's cycle would look like without perturbation. One can see that the perturbations cause small scale (on the order of a day or two) changes in the length of sleep or wake state. However, once allowed to run unperturbed, the system settles back into a regular 24 hour period. Initial wake and sleep times have drifted, but that is to be expected. When run for a longer time, the system will eventually drift back into its original phase with the same wake and sleep initialization times as the subject started with. We know this is an accurate response to these non-drug related, small perturbations(REF). Figures \ref{fig:fig19} and \ref{fig:fig20} show a quantitative measure of the changes in period in the unperturbed vs. the perturbed cycles. Figure \ref{fig:fig21} demonstrates the upholding of limit cycle dynamics in the face of perturbations, and also that the $AD-GABA$ phase plane retains its size and shape but drifts in phase. 

\section{Conclusions}

We have presented a model of the human sleep-wake generating system that identifies specific neuronal regions in the brain with their associated neurotransmitters and projections.  It accounts for the homeostatic influences of sleep-wake, and has the potential to encompass circadian influences as well. It improves upon previous models of sleep-wake cycles in that the oscillations are self-sustaining without the use of external periodic forcing, and it is formulated specially to look at the specific sleep-wake cycle of humans. We use the known neurochemical dynamics of the human sleep-wake cycle to derive equations and form the stability criteria for mathematical phenomena. Since based on limit cycle dynamics, the system runs with the contributions of all the necessary neurotransmitters and with no need of an external driving force. This system gives reasonable predictions for ``natural" sleep bout duration, wake bout duration, REM bout duration and NREM bout duration, as well as estimated natural sleep initiation and wake initiation times. The NREM sleep state is considered the ground state from which all other activity deviates, and REM activity is determined from the difference of the wake and NREM states, as well as the contributions from REM active neurotransmitters. Using experimental data from a mammalian orexin knockout system as well as simply and non-chemically perturbed human sleep-wake cycles, we were able to visualize the perturbed cycle and compare to the experimental results without instability occurring. We can draw a number of conclusions from the resulting model. In the human systems, given forced prolonged waking, sleep duration reacts accordingly to build up appropriate wake pressure. Once the cycle is restored to normal sleep onset and wake onset, the same period for an unperturbed system is maintained. The system exhibits some phase drift with this forcing. Given forced prolonged sleep, wake duration also reacts accordingly to build up appropriate sleep pressure. Once the cycle is restored to normal sleep and wake onset, the periodicity is maintained. The system exhibits some phase drift with this forcing. Given forced wake up during sleep, the subject does not stay awake for long, and falls back into a normal sleep schedule, yielding only some phase drift that is corrected after a few cycles. A subject's REM/NREM cycle goes through three or four oscillations per state (REM or NREM). The first NREM sleep stage is triggered by sleep onset, or when adenosine levels fall below a certain threshold and GABA becomes the dominant neurotransmitter. The last REM cycle is turned off by wake onset, or when adenosine pressure builds back up and GABA decays. Subjects naturally fall asleep and wake up at the intersection points of the GABA and Adenosine thresholds - falling asleep when enough sleep pressure from AD has been accumulated, and waking up when the sleep pressure has decayed once again and when the subject has gone though three to four cycles each of REM and NREM sleep. We can see that given no external forcing (no stimuli, drugs, etc), a ``normal" human's sleep-wake schedule unfolds in a periodic, 24 hour fashion that endures indefinitely unless given some sort of perturbation. \\

It is important to note the implications of the differences between homeostatic and circadian drive, and the importance of both of these drives' effect on sleep. Unlike most recent neurochemically driven models, ours has no need for external forcing to keep it periodically driven. The periodic, limit cycle activity arises naturally from the neurotransmitter dynamics and the equations themselves. This demonstrates an intrinsic, homeostatic force that keeps the human sleep-wake cycle functioning without a change in periodicity, regardless of whether circadian control is included or not. The perturbations in the experimental section demonstrate that even when subjected to strong forcing, the sleep-wake cycle retains its natural 24 hour period, and shifts in phase are the only real artifacts of outside forcing. The next step would be to develop a more complete predictive model of the human sleep wake cycle for use in the scientific and medical fields. \\

\section{Acknowledgements}
This work was supported, in part, though NSF Grant DMS-0636358.  Also, we would like to thank the following for their aid in this project: Dr. Mary Carskadon (Brown University), Dr. Tom Scammell and the Scammell Lab  (BIDMC, Harvard Medical School).

\newpage
\bibliographystyle{plainnat}
\bibliography{project_description_references} 

\end{document}